\documentclass[12pt]{article}
\usepackage{dsfont}
\usepackage{diagbox}
\usepackage{subfigure} 
\usepackage{amsmath, appendix, ulem, amsthm}
\usepackage[nobysame]{amsrefs}
\usepackage{amssymb, color}
\usepackage[margin=1in]{geometry}
\usepackage{mathrsfs}
\usepackage{graphicx}
\usepackage{float}
\usepackage[linesnumbered,ruled]{algorithm2e}
\usepackage{epsf}
\usepackage{titletoc}
\usepackage{cases}
\usepackage{algpseudocode}  
\usepackage{amsmath}  
\graphicspath{{../Figures/}}

\numberwithin{equation}{section}

\theoremstyle{remark}

\renewcommand{\hat}{\widehat}

\newcommand{\f}{\frac}

\newcommand{\beq}{\begin{equation}}
\newcommand{\eeq}{\end{equation}}

\newcommand{\eps}{\varepsilon}

\newcommand{\PP}{\mathcal{P}}

\newcommand{\rd}{\mathrm{d}}

\newcommand{\vpran}[1]{\left(#1\right)}
\usepackage[colorlinks,
            linkcolor=red,
            anchorcolor=red,
            citecolor=red
            ]{hyperref}

\newcommand{\Dx}{\Delta x}

\newcommand{\half}{\frac{1}{2}}
\newcommand{\order}{\mathcal{O}}
\newcommand{\average}[1]{\left\langle#1\right\rangle}

\usepackage{authblk}

\begin{document}
\title{A fully asymptotic preserving decomposed multi-group method for the frequency-dependent radiative transfer equations}
\author[1]{Xiaojiang Zhang\thanks{xjzhang123@sjtu.edu.cn}}
\author[2,3]{Peng Song\thanks{song\_peng@iapcm.ac.cn}}
\author[4]{Yi Shi\thanks{shiyi@sdu.edu.cn}}
\author[5]{Min Tang\thanks{tangmin@sjtu.edu.cn}}
\affil[1,5]{School of Mathematical Sciences, Institute of Natural Sciences, MOE-LSC, Shanghai Jiao Tong University, Shanghai, 200240, P.R. China.}
\affil[2]{ Institute of Applied Physics and Computational Mathematics, Beijing, 100094, P.R. China.} 
\affil[3]{HEDPS, Center for Applied Physics and Technology, College of Engineering, Peking University, Beijing, 100871, P.R. China.}
\affil[4]{School of Mathematics, Shandong University, Jinan, 250100, P.R. China.}

\date{}
\maketitle 

\begin{abstract}
	The opacity of FRTE depends on not only the material temperature but also the frequency, whose values may vary several orders of magnitude for different frequencies. The gray radiation diffusion and frequency-dependent diffusion equations are two simplified models that can approximate the solution to FRTE in the thick opacity regime. The frequency discretization for the two limit models highly affects the numerical accuracy. However, classical frequency discretization for FRTE  considers only the absorbing coefficient. In this paper, we propose a new decomposed multi-group method for frequency discretization that is not only AP in both gray radiation diffusion and frequency-dependent diffusion limits, but also the frequency discretization of the limiting models can be tuned. Based on the decomposed multi-group method, a full AP scheme in frequency, time, and space is proposed. Several numerical examples are used to verify the performance of the proposed scheme.
\end{abstract}


\section{Introduction}
The frequency-dependent radiative transfer equation (FRTE) describes the time evolution of the probability density function of photons that transport and interact with the background material, which plays an important role in astrophysics, inertial confinement fusion (ICF), and high energy physics. Due to the complexity of the FRTE and its wide applications, FRTE has attracted the attention of many physicists and mathematicians.

The numerical simulation of FRTE is challenging due to the following difficulties: 1) The high dimensionality. The photon probability distribution function depends on space, time, moving velocity and photon frequency, which is a 7 dimensional function. For the simplest slab geometry, the dimension is 4, which is hard to find an efficient solver due to the "curse of dimensionality"; 2) Multiscale parameters.
Absorbing and scattering coefficients are two important material properties that influence the photons' movements. The absorbing  coefficient $\sigma_a$ (scattering coefficient $\sigma_s$) gives the probability that photons will be absorbed (scattered) by the background material. Both $\sigma_a$ and $\sigma_s$ may depend on the material type, material temperature and the photon frequency. When $\sigma_a+\sigma_s$ is large, the mean free path (the average distance between two successive scattering or absorption) of the photon is small, the material is referred to as "optical thick", while when $\sigma_a+\sigma_s$ events is at $O(1)$, the material is "optical thin". For photons with different frequencies, the values of $\sigma_a$ and $\sigma_s$ may vary several orders of magnitude, which induces the multiscale property of the system. To guarantee both accuracy and stability, classical schemes require very fine space and time steps, which leads to high computational costs. 3) The strong nonlinearity. The material opacity and the Planck function which gives the equilibrium radiation frequency distribution depend nonlinearly on the material temperature. If one needs an efficient solver that allows for large time steps, some terms have to be treated implicitly and thus nonlinear iterations are required. However, the nonlinear iterations for high dimensional equations are extremely expensive.

One popular strategy for the FRTE simulation is the stochastic method \cite{steinberg2022multi,densmore2012hybrid} based on the Fleck-Cummings' implicit Monte Carlo (IMC) method \cite{FC71}. However, it requires a large number of sampling particles and suffers from unavoidable statistical fluctuations. On the other hand, there are deterministic solvers for the radiative transfer equation (RTE), which is a simplified version of FRTE that ignores the frequency dependence. Deterministic solvers include Discontinuous Galerkin method \cite{lcs,bailey2008piecewise}, photon free method \cite{chang2007incorporation,chang2007deterministic}, diffusion-synthetic acceleration method (DSA) \cite{H1993,Larsen1988,ALarsen2002} etc. and some of them have been extended to FRTE.


In order to deal with the multiscale parameters in FRTE, two strategies are proposed recently. One is designing asymptotic preserving (AP) schemes. The AP scheme provides a general framework of solving the difficulty coming from multiscale parameters. A scheme is AP when the asymptotic limit of the discretization  becomes  a stable  solver for the limit model as the scaling parameter in the microscopic equation goes to zero \cite{jin20101,hu2017asymptotic,2018Tang}. In last decade, research on the AP schemes has been developed rapidly for various applications and we refer interesting readers to a recent review paper \cite{jin_2022}.  The other approach is the moment-based acceleration scheme which is called High-Order/Low-Order (HOLO) algorithm \cite{chacon2017multiscale,yee2017stable,park2019multigroup}. In the HOLO algorithm, a low-order system which consists of the first two  moments of the RTE and the material temperature equation is solved. Then the RTE is solved by determining the Planckian emission source term using the temperature obtained in the lower-order system. The benefit of HOLO algorithms is nonlinear elimination, but their long-term accuracy and nonlinear stability require discrete consistency between HO and LO formulations \cite{chacon2017multiscale}. Both approaches are well developed for RTE simulations.

Designing AP schemes for the nonlinear radiation hydrodynamic system is a popular topic recently. For the gray radiative transfer equations (GRTE) in the gray radiative diffusion limit, several AP schemes can be found in the literature. An AP scheme is constructed in \cite{2001Klar} by decomposing the distribution function into the equilibrium and non-equilibrium parts; in \cite{2020park}, the authors developed an AP-HOLO algorithm based on  the linear discontinuous Galerkin HO transport solver and the corner-balance LO solver; in \cite{sun2015asymptotic}, an unified gas kinetics scheme (AP-UGKS) was proposed. However, in order to use an opacity independent time step, a linearized  iterative  solver involving the radiation equation has to be employed for the time-implicit part.
In \cite{tang2021accurate}, an AP accurate front capturing scheme was proposed that allows both space and time steps being independent of the opacity and only requires a scalar Newton's solver. 
Later on, in order to get AP schemes for more complex applications, UGKS has been extended to FRTE \cite{sunj} and the radiation magnetohydrodynamics (RMHD) system \cite{jin2022spatial,sun2020multiscale}; The AP linear-discontinuous spatial differencing scheme has been extended to FRTE in \cite{MWTS}; The idea of intensity decomposition used in \cite{tang2021accurate} has been extended to the RMHD system in \cite{jin2022spatial}.
However, none of the above mentioned work considered about the improvement of the frequency discretization.

\cite{terh} separates particles into several groups. It is assumed that particles in the same group have the same frequency and the same absorption and scattering coefficients. However, due to the limitation of computational resources, the number of groups is usually not enough to resolve all variances in the frequency distribution. 
On the other hand, when the opacity is thick, there are two approximating models, one is the gray radiation diffusion model where the particles are at the equilibrium Planck distribution in frequency; the other is the frequency-dependent diffusion limit (FDDL), where particles with different frequencies diffuse at different speeds. In order to use a small number of groups to provide a correct solution behavior,  different frequency discretizations of the limit models are proposed in the literature, \textcolor{black}{for example,  the piece-constant approximation \cite{chang2007deterministic,chang2007incorporation}, the Rosseland mean approximation \cite{MWTS,yee2017stable}, the Planck approximation \cite{sunj,sunj18,2017An}.} However, it is not clear how to construct a group discretization for the FRTE whose diffusion limits can coincide with the group discretization of the limiting models.    

In this paper, by decomposing the intensity into three parts, we obtain a new decomposed multi-group method whose asymptotic limit, when the mean free path goes to zero, can be consistent with the required
frequency discretization of both limiting models. To the best of our knowledge, the decomposed multi-group method is the first work that gives an AP discretization in the frequency domain.

Based on the new decomposed multi-group method, a full 
 AP scheme in frequency, space, and time is proposed. The scheme accuracy is guaranteed in both optically thin and thick regimes when large spatial and time steps are used.
Moreover, nonlinear iterations are employed to solve a system with only macroscopic quantities. Then linear transport equations with decoupled direction and energy groups are solved. Thus the computational cost is similar to HOLO \cite{park2019multigroup}.


The organization of this paper is as follows. Section 2 shows two different asymptotic limits of the FRTE, one is the gray radiative diffusion limit, and the other is the FDDL. The classical multi-group method and our decomposed multi-group method are compared in the frequency-dependent optical thick regime in Sections 3 and 4. The semi-discretization that is AP in frequency and time is proposed in section 5, and the fully discretized scheme that is AP in frequency, time, and space is displayed in section 6. Their AP properties in both gray radiation and frequency-dependent diffusion limits are shown by asymptotic analysis. In section 7, some numerical tests are given to illustrate the stability and accuracy of our proposed scheme. Finally, the paper is concluded with some discussions in Section 8.

\section{The FRTE and its limit models}\label{frte}
\subsection{The model and nondimensionalization}
The FRTE under consideration is as in \cite{2017An}, that writes:
\begin{equation}\label{eqn:001}  
\left\{
\begin{aligned}
&\frac{1}{c}\partial_t I+\Omega \cdot \nabla I  = \sigma_a\left( B(\nu,T)-I\right)+\sigma_s\left( \rho-I\right), \qquad \rho=\f{1}{4\pi}\int_{4\pi}I \text{d}\Omega := \average{I},\\
& C_v\partial_t T\equiv  \partial_t U  =  \int\limits_{4\pi}\int\limits_0^{\infty } \sigma_a\left( I-B(\nu,T)\right)\text{d}\nu \text{d}\Omega,\ \  \\
\end{aligned}
\right.
\end{equation} 
where $x$ is the spatial variable; $\Omega$ is the angular variable; $\nu\in(0,+\infty)$  is the frequency variable; $I(t,x,\Omega,\nu)$ is the radiation intensity; $T(t,x)$ is the material temperature; $\sigma_a(x,\nu,T)$ and $\sigma_s(x,\nu,T)$ are the absorbing and scattering coefficient respectively; $c$ and $C_v$ represent respectively the light speed and the specific heat capacity; the Planck function $B(\nu,T)$ is defined by
\begin{equation*}
	B(\nu,T)=\frac{2h\nu^3}{c^2}\frac{1}{e^{h\nu/kT}-1},
\end{equation*}
with  Boltzmann's constant $k$ and Planck’s constant $h$. The Planck function satisfies
$$
 \int\limits_0^{\infty }4\pi B \ \text{d}\nu=a_r cT^4, \qquad  \int\limits_0^{\infty }4\pi\frac{\partial B}{\partial T} \ \text{d}\nu=4a_r cT^3,
$$
 with the radiation constant $a_r$. The FRTE system \eqref{eqn:001} is composed of two equations, the first one describes the time evolution of radiation intensity and the second one gives the time evolution of background material temperature. If we integrate the first equation in system \eqref{eqn:001} with respect to the frequency $\nu$ from $0$ to $+\infty$ and assume that the absorbing and scattering coeﬀicients are independent of the frequency, the GRTE can be obtained. GRTE is a good approximation to FRTE when the media is assumed to be gray, i.e. its properties are independent of the radiation spectrum \cite{traugott1968grey}. However, in reality, media opacity varies a lot for different temperatures, densities and particle frequencies. FRTE has to be considered in some situations  \cite{vaytet2011numerical,vaytet2013influence}.

We then consider the dimensionless form of \eqref{eqn:001}. Similar as in \cite{LREJ,HFJM}, we consider the following nondimensionalization:
$$
x=\hat{x}\ell_\infty,\quad t=\hat{t} t_\infty,\quad \nu=\hat{\nu}\f{kT_\infty}{h},\quad \sigma_a=\hat{\sigma}_a/\ell_\infty, \quad\sigma_s=\hat{\sigma}_s/\ell_\infty,
$$
$$
I=\hat{I}\f{a_rchT_\infty^3}{k},\quad B=\hat{B}\f{a_rchT_\infty^3}{k},\quad \rho=\hat{\rho}\f{a_rchT_\infty^3}{k},\quad U=\hat{U}\f{\ell_\infty^2}{t_\infty^2},
$$
where  variables with a hat denote nondimensional quantity and variables with $\infty- $subscript are the characteristic value with unit. More precisely, $\ell_\infty$, $t_\infty$ and $T_\infty$ are respectively the reference length, time and temperature, and $a_r$ is the radiation constant. 

Then the full dimensionless FRTE  becomes 
\begin{equation}\label{eqn:001-lessdim}  
\left\{
\begin{aligned}
&\frac{1}{\mathcal{C}}\partial_{\hat{t}} \hat{I}+\Omega \cdot \nabla_{\hat{x}} \hat{I}  = \hat{\sigma}_a\mathscr{L}_a\left( \hat{B}-\hat{I}\right)+\hat{\sigma}_s\mathscr{L}_s\left( \hat{\rho}-\hat{I}\right),\\
& \hat{C}_v \partial_{\hat{t}} \hat{T}\equiv \partial_{\hat{t}} \hat{U}=  \mathcal{C}\mathcal{P}_0\int\limits_{4\pi}\int\limits_0^{\infty } \hat{\sigma}_a\mathscr{L}_a\left( \hat{I}-\hat{B}\right)\text{d}\hat{\nu} \text{d}\Omega,\ \  \\
\end{aligned}
\right.
\end{equation} 
where $\mathcal{C}=\f{ct_\infty}{\ell_\infty}$, $\mathcal{P}_0=\f{a_rt_\infty^2T_\infty^4}{\ell^2_\infty}$, $\hat{C}_v$ is a positive constant and $\hat{\sigma}_a\mathscr{L}_a$, $\hat{\sigma}_s\mathscr{L}_s$ are respectively $\sigma_a$ and $\sigma_s$. According to \cite{LMH}, $\mathcal{P}_0$  is a nondimensional parameter that measures radiation effects on the
	matter, which is always considered as $\mathcal{P}_0=\order(1)$ when a moderate amount of radiation in the matter is considered. $\mathcal{C}=\order(1/\eps)$ for nonrelativistic case,  which indicates that the speed of light travels much faster than the material velocity when $\eps\ll 1$. Thus, in this paper, we fix the values of $\mathcal{P}_0$ and $\mathcal{C}$ as \begin{equation}\mathcal{P}_0=1,\qquad \mathcal{C}=1/\eps.\label{eq:p0c}\end{equation}
Depending on the scales of $\mathscr{L}_a$ and $\mathscr{L}_s$, the solution of system \eqref{eqn:001} tends to solutions of different limit systems in the diffusive regime \cite{GPG}. How the different scales of $\mathscr{L}_a$ and $\mathscr{L}_s$ lead to different limits will be discussed in the subsequent part. For simplicity,  the hat in \eqref{eqn:001-lessdim} will be dropped later on.

\subsection{The diffusion limits of the FRTE}\label{d-frte}
 We focus on two different regimes for the FRTE, one is the gray radiation diffusion regime and the other is the frequency dependent diffusion regime, whose derivations are similar as in \cite{GPG,LMH}.  The gray radiation diffusion limit is valid when collisions are so frequent that the frequency distribution depends predominantly on the local temperature. In the FDDL, the determination of the frequency distribution becomes complicated, it depends on both the local temperature and the state the radiation field \cite{Holland1969,manuel2001non}. Both limiting models are widely used in the simulations of RMHD problems \cite{CHENG2020109724,CD20,jin2022spatial}. 
 For the convenience of readers, the details of the derivation are in  Appendix \ref{apen01}.
 
 Using \eqref{eq:p0c}, the FRTE \eqref{eqn:001} becomes:
\begin{subequations}\label{eqn:001-eq}  
\begin{numcases}{}
\eps\partial_t I+\Omega \cdot \nabla I  =\mathscr{L}_a\sigma_a\left( B(\nu,T)-I\right)+\mathscr{L}_s\sigma_s\left( \rho(t,x,\nu)-I\right), \label{eqn:001-eq01}\\
 \eps C_v\partial_t T  =  \int\limits_{4\pi}\int\limits_0^{\infty }\mathscr{L}_a\sigma_a\left( I-B(\nu,T)\right)\text{d}\nu \text{d}\Omega.\label{eqn:001-eq02}
\end{numcases}
\end{subequations}

\begin{itemize}
\item \textbf{The gray radiation diffusion regime.} In this regime,  the radiation intensity $I$ adapts to the material temperature. The following three different scalings yield the gray radiation diffusion equation:
\begin{subequations}\label{eqn:scaling}  
	\begin{numcases}{}
\mathscr{L}_a = 1/\eps,\quad\mathscr{L}_s = \eps, \label{eqn:sca01}\\
\mathscr{L}_a = 1,\quad\mathscr{L}_s =1/ \eps,\label{eqn:sca02}\\
\mathscr{L}_a = 1/\eps,\quad\mathscr{L}_s =1/ \eps.\label{eqn:sca03}
	\end{numcases}
\end{subequations}
When $\eps\to 0$ in (\ref{eqn:001-eq}), the material temperature $T$ can be approximated by the solution of the following gray radiation diffusion equation
 \begin{equation}\label{eq:limit_T}
\partial_t (T^{(0)})^4+C_v \partial_t T^{(0)}=\nabla \cdot \vpran{\frac{D_d}{\sigma  }\nabla (T^{(0)})^4}\,,
\end{equation}
with $D_d=\f{1}{3}I_d$  (where $I_d$ denotes the 3 by 3 identity matrix) and
the mean opacity $\sigma$ being given by
\begin{equation}\label{eq:defsigma}
\frac{1}{\sigma(x, T)} \equiv \frac{\int_{0}^{\infty} \frac{4\pi\mathcal{C}}{\mathscr{L}_a\sigma_a(x, \nu, T)+\mathscr{L}_s\sigma_s(x, \nu, T)} \frac{\partial B^{(0)}(\nu, T)}{\partial T} \text{d} \nu}{\int_{0}^{\infty}4\pi \frac{\partial B^{(0)}(\nu, T)}{\partial T} \text{d} \nu}=\frac{\int_{0}^{\infty} \frac{4\pi\mathcal{C}}{\mathscr{L}_a\sigma_a(x, \nu, T)+\mathscr{L}_s\sigma_s(x, \nu, T)} \frac{\partial B^{(0)}(\nu, T)}{\partial T} \text{d} \nu}{4  (T^{(0)})^{3}}.
\end{equation}
The radiation intensity $I$ can then be approximated by \begin{equation}\label{appi}
B(\nu,T^{(0)})-\f{1}{\mathscr{L}_a\sigma_a+\mathscr{L}_s\sigma_s}\f{\partial B(\nu,T^{(0)})}{\partial T}\Omega\cdot\nabla T^{(0)}.
\end{equation}

\item \textbf{The frequency dependent diffusion regime.} In this regime, the radiation intensity $I$ tends to be independent of $\Omega$, but the frequency distribution can no longer be determined by the material temperature $T$.
Similar as in \cite{GPG}, we consider the following scaling
$$
\mathscr{L}_a = \eps,\quad\mathscr{L}_s = 1/\eps.
$$
 In this regime, the FRTE \eqref{eqn:001-eq} writes:
\begin{subequations}\label{eqn:001-noneq}  
\begin{numcases}{}
\eps\partial_t I+\Omega \cdot \nabla I  = \eps\sigma_a\left( B(\nu,T)-I\right)+\f{\sigma_s}{\eps}\left( \rho(t,x,\nu)-I\right),\label{eqn:001-noneq01}\\
 C_v\partial_t T  =  \int\limits_{4\pi}\int\limits_0^{\infty }\sigma_a\left( I-B(\nu,T)\right)\text{d}\nu \text{d}\Omega.\label{eqn:001-noneq02}
\end{numcases}
\end{subequations}
When $\eps\to 0$ in (\ref{eqn:001-noneq}), the solution $I(t,x,\Omega,\nu)\approx \rho^{(0)}(t,x,\nu)$ where $\rho^{(0)}$ and $T$ satisfies the following frequency dependent diffusion system:
\begin{equation}\label{nonequ-diff}  
\left\{
\begin{aligned}
&\partial_t \rho^{(0)}- \nabla \cdot \vpran{\frac{D_d}{\sigma_{s}  }\nabla \rho^{(0)}}=\sigma_a(B^{(0)}-\rho^{(0)}),\\
&  C_v\partial_t T^{(0)}  =\int_{0}^{\infty}4\pi\sigma_a(\rho^{(0)}-B^{(0)}) \text{d} \nu.\ \  \\
\end{aligned}
\right.
\end{equation} 
\end{itemize}


\section{Extension of the classical multi-group method and the diffusive limits}\label{apen02}
\subsection{The classical multi-group method and its extension}
We first consider the frequency discretization of the system \eqref{eqn:001-lessdim}. In the multi-group method, the continuous frequency space $(0,\infty)$ is divided into $G$ groups, and each frequency interval is denoted by $(\nu_{g-1/2},\nu_{g+1/2})$ ($g=1,\cdots, G$), with $\nu_{\frac{1}{2}}=0$, $\nu_{G+\frac{1}{2}}=+\infty$.  By integrating the first equation in system (\ref{eqn:001-lessdim}) over the frequency interval $(\nu_{g-1/2},\nu_{g+1/2})$ yields
\begin{equation}\label{eqn:002}  
\int\limits_{\nu_{g-1/2}}^{\nu_{g+1/2}}\f{1}{\mathcal{C}}\partial_t I+\Omega \cdot \nabla I \text{d}\nu = \int\limits_{\nu_{g-1/2}}^{\nu_{g+1/2}}\mathscr{L}_a\sigma_a\left( B(\nu,T)-I\right)+\mathscr{L}_s\sigma_s\left( \rho-I\right)\text{d}\nu.\\
\end{equation} 
Then the multi-group discretization of FRTE becomes:
\begin{subequations}\label{eqn:003re}  
	\begin{numcases}{}
	\f{1}{\mathcal{C}}\partial_t I_g+\Omega \cdot \nabla I_g = \int\limits_{\nu_{g-1/2}}^{\nu_{g+1/2}}\mathscr{L}_a\sigma_a\left( B-I\right)-\mathscr{L}_s\sigma_s\left(  \rho-I\right)\text{d}\nu,\quad (g=1,\cdots,G),  \label{mgre1}\\
	C_v\partial_t T  = \mathcal{C}\PP_0 \int\limits_{4\pi}\sum\limits_{g'=1}^{G }\int\limits_{\nu_{g-1/2}}^{\nu_{g+1/2}} \left[\mathscr{L}_a\sigma_{a} \left(  I-  B\right)\right]\text{d}\nu\text{d}\Omega,\label{mgre2}
	\end{numcases}
\end{subequations}
in which the radiation intensity $I(t,x,\Omega,\nu)$ and Planck function $B(\nu,T)$ in different groups  are given by
\begin{equation}\label{defbg}
I_{g}=\int\limits_{\nu_{g-\frac{1}{2}}}^{\nu_{g+\frac{1}{2}}} I(t, x, \Omega, \nu) \text{d} \nu, \qquad B_{g}=\int\limits_{\nu_{g-\frac{1}{2}}}^{\nu_{g+\frac{1}{2}}} B(\nu,T) \text{d} \nu.
\end{equation}

After frequency discretization, only $I_g$ in each group and the Planck function $B$ are known, while the full frequency distribution of $I$ is not known. Thus the main difficulty is how to approximate the two terms
\begin{equation}\label{eq:sigmaasiamgsdecomposed}
\int\limits_{\nu_{g-1/2}}^{\nu_{g+1/2}}\sigma_a\left( B-I\right)\text{d}\nu\quad\mbox{and}\quad \int\limits_{\nu_{g-1/2}}^{\nu_{g+1/2}}\sigma_s\left(  \rho-I\right)\text{d}\nu.
\end{equation}
In order to approximate the two terms on the right hand side of \eqref{eqn:002} by $I_g$, $B_g$ and $T$, some approximations can be found in the literature. We review some classical multi-group frequency discretization for FRTE in the subsequent part. Most works in the literature consider only the case when $\sigma_s=0$. In this subsection, we extend similar idea of approximating the term $\int\limits_{\nu_{g-1/2}}^{\nu_{g+1/2}}\mathscr{L}_a\sigma_a\left( B(\nu,T)-I\right)$ in \cite{yee2017stable,MWTS} to the case when $\sigma_s\neq 0$.

In the classical multi-group frequency discretization in \cite{terh}, by defining $\sigma_{a,g}$ and  $\sigma_{s,g}$ as follows
\begin{equation}\label{eqapp}
\sigma_{a, g}=\frac{\int_{\nu_{g-\frac{1}{2}}}^{\nu_{g+\frac{1}{2}}} \sigma_{a}(B-I) \mathrm{d} \nu}{B_{g}-I_{g}}, \quad \sigma_{s, g}=\frac{\int_{\nu_{g-\frac{1}{2}}}^{\nu_{g+\frac{1}{2}}} \sigma_{s}(\rho-I) \mathrm{d} \nu}{\rho_{g}-I_{g}},
\end{equation}the multi-group frequency discretization for FDDL writes
\begin{subequations}\label{eqn:003re-tr}  
	\begin{numcases}{}
		\f{1}{\mathcal{C}}\partial_t I_g+\Omega \cdot \nabla I_g = \mathscr{L}_a\sigma_{a,g}\left( B_g-I_g\right)+\mathscr{L}_s\sigma_{s,g}\left(  \rho_g-I_g\right),\qquad (g=1,\cdots,G),  \label{mgre1-tr}\\
		C_v\partial_t T  = \mathcal{C}\PP_0 \int\limits_{4\pi}\sum\limits_{g'=1}^{G } \left[\mathscr{L}_a\sigma_{a,g} \left(  I_g-  B_g\right)\right]\text{d}\Omega.\label{mgre2-tr}
	\end{numcases}
\end{subequations}
Here $\sigma_{a,g}$ and $\sigma_{s,g}$ are respectively the mean absorption and scattering coefficients that have to be approximated by $I_g$, $B_g$ and $T$. In the diffusive regime,  different approximations yield different coefficients in the limit equations. Two classical approximations are discussed in the subsequent part.

The most straight forward way is to use piece-wise constant approximation and let
$B(\nu)\approx B_g$ and $\rho(\nu)\approx \rho_g$, $I(\nu)\approx I_g$ when $\nu\in(\nu_{g-\frac{1}{2}},\nu_{g+\frac{1}{2}})$. Then \begin{equation}\sigma_{a,g}\approx \f{1}{\nu_{g+1/2}-\nu_{g-1/2}}\int_{\nu_{g-\frac{1}{2}}}^{\nu_{g+\frac{1}{2}}} \sigma_{a}\mathrm{d} \nu,\qquad \sigma_{s,g}\approx \f{1}{\nu_{g+1/2}-\nu_{g-1/2}}\int_{\nu_{g-\frac{1}{2}}}^{\nu_{g+\frac{1}{2}}} \sigma_{s}\mathrm{d} \nu.\label{eq:piecewisesection}\end{equation}

In \cite{yee2017stable,MWTS} and some related work later on, the authors considered the Rosseland mean for which $\sigma_{a,g}$ are determined from the approximation in \eqref{appi} which yields

$$
B-I\approx\f{1}{\mathscr{L}_a\sigma_a+\mathscr{L}_s\sigma_s}\f{\partial B}{\partial T}(T)\Omega\cdot\nabla T,
$$
and the Rosseland mean is
\begin{equation}\label{eqn::ra}
\sigma_{a, g}=\frac{\int_{\nu_{g-\frac{1}{2}}}^{\nu_{g+\frac{1}{2}}}\f{\sigma_a}{\mathscr{L}_a\sigma_a+\mathscr{L}_s\sigma_s} \frac{\partial B}{\partial T}(T) \mathrm{~d} \nu}{\int_{\nu_{g-\frac{1}{2}}}^{\nu_{g+\frac{1}{2}}} \f{1}{\mathscr{L}_a\sigma_a+\mathscr{L}_s\sigma_s} \frac{\partial B}{\partial T} (T)\mathrm{~d} \nu}.
\end{equation}
Extending similar idea to the approximation of $\sigma_{s,g}$, the most straight forward way is to use the following approximation
$$
\rho-I\approx-\f{1}{\mathscr{L}_a\sigma_a+\mathscr{L}_s\sigma_s}\Omega\cdot\nabla\rho.
$$
If the intensity $I$ has Planck distribution in the frequency domain, i.e.
$\rho\approx I\approx B(\nu,T_r)$,
where the radiation temperature $T_r$ is given by
\begin{equation}\label{radiat}
T_r^4=\int_{4\pi}\int_{0}^{+\infty}I\mathrm{d} \Omega\mathrm{d} \nu=4\pi\int_0^\infty\rho\mathrm{d}\nu=4\pi\sum\limits_{g'=1}^{G } \rho_g,
\end{equation}
the scattering coefficient $\sigma_{s,g}$ can then be approximated by 
\begin{equation}\label{ros-s}
\begin{aligned}
\sigma_{s, g}&=\frac{\int_{\nu_{g-\frac{1}{2}}}^{\nu_{g+\frac{1}{2}}} \sigma_{s}(\rho-I) \mathrm{d} \nu}{\int_{\nu_{g-\frac{1}{2}}}^{\nu_{g+\frac{1}{2}}} (\rho-I) \mathrm{~d} \nu}\approx\frac{\int_{\nu_{g-\frac{1}{2}}}^{\nu_{g+\frac{1}{2}}}\f{\sigma_s}{\mathscr{L}_a\sigma_a+\mathscr{L}_s\sigma_s} \Omega \cdot\nabla \rho \mathrm{d} \nu}{\int_{\nu_{g-\frac{1}{2}}}^{\nu_{g+\frac{1}{2}}} \f{1}{\mathscr{L}_a\sigma_a+\mathscr{L}_s\sigma_s}\Omega \cdot\nabla \rho \mathrm{~d} \nu}\\
&=\frac{\int_{\nu_{g-\frac{1}{2}}}^{\nu_{g+\frac{1}{2}}}\f{\sigma_s}{\mathscr{L}_a\sigma_a+\mathscr{L}_s\sigma_s} \frac{\partial B}{\partial T}(T_r) \Omega \cdot\nabla T_r \mathrm{d} \nu}{\int_{\nu_{g-\frac{1}{2}}}^{\nu_{g+\frac{1}{2}}} \f{1}{\mathscr{L}_a\sigma_a+\mathscr{L}_s\sigma_s}\frac{\partial B}{\partial T}(T_r)\Omega \cdot\nabla T_r \mathrm{~d} \nu}
=\frac{\int_{\nu_{g-\frac{1}{2}}}^{\nu_{g+\frac{1}{2}}}\f{\sigma_s}{\mathscr{L}_a\sigma_a+\mathscr{L}_s\sigma_s} \frac{\partial B}{\partial T}(T_r) \mathrm{~d} \nu}{\int_{\nu_{g-\frac{1}{2}}}^{\nu_{a}} \f{1}{\mathscr{L}_a\sigma_a+\mathscr{L}_s\sigma_s} \frac{\partial B}{\partial T}(T_r) \mathrm{~d} \nu}.
\end{aligned}
\end{equation}
It is easy to see that the Rosseland mean is an approximation to the average of $\sigma_a$, $\sigma_s$ inside the interval $(\nu_{g+\frac{1}{2}},\nu_{g+\frac{1}{2}})$ as well.

\subsection{The  diffusion limit of the classical multi-group method}
\subsubsection{The gray radiation diffusion regime }
When $\mathscr{L}_a=1/\eps$, $\mathscr{L}_s=1/\eps$, the multi-group frequency discretization in \eqref{eqn:003re-tr} becomes
	\begin{subequations}\label{eqn:003re-equi-s0}  
		\begin{numcases}{}
			\eps\partial_t I_g+\Omega \cdot \nabla I_g = \f{\sigma_{a,g}}{\eps}\left(  B_g- I_g\right)+\f{\sigma_{s,g}}{\eps}\left(  \rho_g- I_g\right),\qquad (g=1,\cdots,G),   \label{mgre1-equi-s0}\\
		\eps C_v\partial_t T  = 4\pi\sum\limits_{g'=1}^{G } \left[\f{\sigma_{a,g'}}{\eps} \left(  \rho_{g'}-  B_{g'}\right)\right].\label{mgre2-equi-s0}
		\end{numcases}
	\end{subequations}
	We derive its gray radiation diffusion limit equation in the subsequent part.
	By using Chapman-Enskog expansion in \eqref{mgre1-equi-s0}, one has
	\begin{equation}\label{eqapp1-appen-s0}
\begin{aligned}
	I_g &=\frac{ \sigma_{a,g}}{ \sigma_{a,g} +\sigma_{s,g}} B_g+\frac{\sigma_{s,g}}{  \sigma_{a,g} +\sigma_{s,g}} \rho_g-\frac{\varepsilon}{ \sigma_{a,g} +\sigma_{s,g}}\left[\eps\partial_{t} I_g+ \Omega \cdot \nabla I_g\right] \\
&=\frac{ \sigma_{a,g}}{ \sigma_{a,g} +\sigma_{s,g}} B_g+\frac{\sigma_{s,g}}{  \sigma_{a,g} +\sigma_{s,g}} \rho_g-\frac{\varepsilon}{\sigma_{a,g} +\sigma_{s,g}} \Omega \cdot \nabla \left(\frac{ \sigma_{a,g}}{ \sigma_{a,g} +\sigma_{s,g}} B_g+\frac{\sigma_{s,g}}{  \sigma_{a,g} +\sigma_{s,g}} \rho_g\right)+\mathcal{O}\left(\varepsilon^{2}\right).
\end{aligned}
\end{equation}
	Taking the integral with respect to $\Omega$ of equation \eqref{mgre1-equi-s0}, we can obtain
	\[
	\rho_g =  B_g + \order(\eps)\,,
	\]
	which implies that $\rho_g\approx B_g$.  Taking the integral with respect to $\Omega$ of the equation \eqref{mgre1-equi-s0}, adding up all groups and  \eqref{mgre2-equi-s0}, one can get
	\begin{equation} \label{diff-mgre-s0}
	4\pi\partial_t\left(\sum_{g=1}^G \rho_g\right)  +\f{4\pi}{\eps} \sum_{g=1}^G\nabla \cdot \average{\Omega I_g} + C_v \partial_t T = 0\,.
	\end{equation}
	Plugging \eqref{eqapp1-appen-s0} into the equation \eqref{diff-mgre-s0}, using the condition $\rho_g =  B_g + \order(\eps)$,
	letting $\eps\to 0$, then \eqref{diff-mgre-s0} reduces to
	\begin{equation} \label{diff-mg1re-s0}
	\partial_t\left(\sum_{g=1}^G 4\pi B_g\right) + C_v\partial_t T=  \sum_{g=1}^G\nabla \cdot\left(  \frac{4\pi}{3(\sigma_{a,g}+\sigma_{s,g})} \nabla B_g\right)  \,,
	\end{equation}
	which implies that
		\begin{equation} \label{diff-mg1re-s01}
	\partial_t\left(\sum_{g=1}^G 4\pi B_g\right) + C_v\partial_t T=  \sum_{g=1}^G\nabla \cdot\left(  \frac{4\pi}{3(\sigma_{a,g}+\sigma_{s,g})}\frac{\partial B_g}{\partial T}\frac{1}{4 T^3} \nabla T^4\right)  \,,
	\end{equation}
	When the group discretization of the Planck equilibrium $B$ satisfies the following relations 
\begin{equation}\label{Bgconditions}
4\pi\sum_{g=1}^G  B_g = \int\limits_0^{\infty }4\pi B \ \text{d}\nu=T^4,
\end{equation}
	  \eqref{diff-mg1re-s01} then reduces to
	\begin{equation}
	\partial_t T^4  + C_v\partial_t T=  \nabla \cdot \left(  \frac{1}{3\widehat\sigma_t} \nabla T^4 \right) \,,
	\end{equation}
When one approximates $\sigma_{a,g}$ and $\sigma_{s,g}$ as in \eqref{eq:piecewisesection}, the diffusion coefficient of the above equation is determined by
\begin{equation}\label{sigmatconstant}
\f{1}{\hat\sigma^c_t}= \f{1}{4T^3} \left(\sum\limits_{g=1}^{G }  \frac{4\pi}{\f{1}{\nu_{g+1/2}-\nu_{g-1/2}}\int_{\nu_{g-1/2}}^{\nu_{g+1/2}} {\sigma_a+\sigma_s}  \text{d} \nu} \f{\partial B_g}{\partial T}\right).
\end{equation}
On the other hand,	when one uses the Rosseland mean as in \eqref{eqn::ra} and \eqref{ros-s}, one has
$$
\f{1}{\sigma_{a,g}+\sigma_{s,g}}=\frac{\int_{\nu_{g-\frac{1}{2}}}^{\nu_{g+\frac{1}{2}}}\f{1}{\sigma_a+\sigma_s} \frac{\partial B}{\partial T}(T) \mathrm{~d} \nu}{\int_{\nu_{g-\frac{1}{2}}}^{\nu_{g+\frac{1}{2}}}  \frac{\partial B}{\partial T} (T)\mathrm{~d} \nu},
$$
which implies 
	\begin{equation}\label{eq:sigmat1}
	\f{1}{\widehat\sigma^r_t}=\f{1}{4T^3} \left(\sum\limits_{g=1}^{G }\int\limits_{\nu_{g-1/2}}^{\nu_{g+1/2}}  \frac{4\pi}{\sigma_a+\sigma_s} \f{\partial B}{\partial T}\text{d} \nu\right).
	\end{equation}
Both $\hat\sigma_t^c$ and $\hat\sigma_t^r$ provide reasonable approximations to $\sigma$ in \eqref{eq:defsigma} when the $\sigma_a$, $\sigma_s$ are at $O(1)$ and when there are enough groups. 
	However, we will see in the next subsection that $\hat \sigma^c_t$ determined by piece-wise constant approximation as in \eqref{sigmatconstant} can not provide a good approximation to the the diffusion coefficient as in \eqref{eq:defsigma}, while the
	Rosseland mean $\hat \sigma^r_t$ can. Moreover,
	When $\sigma_s=0$, $\f{1}{\widehat\sigma^r_t}$ becomes
		$$
	\f{1}{\widehat\sigma^r_a}=\f{1}{4T^3} \left(\sum\limits_{g=1}^{G }\int\limits_{\nu_{g-1/2}}^{\nu_{g+1/2}}  \frac{4\pi}{\sigma_a} \f{\partial B}{\partial T}\text{d} \nu\right).
	$$
which is the same as in \cite{MWTS}.

\subsubsection{The frequency dependent  diffusion regime}
When $\mathscr{L}_a=\eps$, $\mathscr{L}_s=1/\eps$, the multi-group frequency discretization in \eqref{eqn:003re-tr} becomes
\begin{subequations}\label{eqn:003re-equi-s0f}  
	\begin{numcases}{}
	\eps\partial_t I_g+\Omega \cdot \nabla I_g = \eps\sigma_{a,g}\left(  B_g- I_g\right)+\f{\sigma_{s,g}}{\eps}\left(  \rho_g- I_g\right),\qquad (g=1,\cdots,G),   \label{mgre1-equi-s0f}\\
	\eps C_v\partial_t T  = 4\pi\sum\limits_{g'=1}^{G } \left[\eps\sigma_{a,g'}\left(  \rho_{g'}-  B_{g'}\right)\right].\label{mgre2-equi-s0f}
	\end{numcases}
\end{subequations}
We derive its  radiation diffusion limit equation in the subsequent part.
By using Chapman-Enskog expansion in \eqref{mgre1-equi-s0f}, one has
\begin{equation}\label{eqapp1-appen-s0f}
\begin{aligned}
I_g &=\frac{ \eps^2\sigma_{a,g}}{\eps^2 \sigma_{a,g} +\sigma_{s,g}} B_g+\frac{\sigma_{s,g}}{  \eps^2 \sigma_{a,g} +\sigma_{s,g}} \rho_g-\frac{\varepsilon}{ \eps^2 \sigma_{a,g} +\sigma_{s,g}}\left[\eps\partial_{t} I_g+ \Omega \cdot \nabla I_g\right] \\
&=\frac{\sigma_{s,g}}{ \eps^2 \sigma_{a,g} +\sigma_{s,g}} \rho_g-\frac{\varepsilon}{\eps^2\sigma_{a,g} +\sigma_{s,g}} \Omega \cdot \nabla \left(\frac{\sigma_{s,g}}{ \eps^2  \sigma_{a,g} +\sigma_{s,g}} \rho_g\right)+\mathcal{O}\left(\varepsilon^{2}\right)\\
&=\rho_g-\frac{\varepsilon}{ \sigma_{s,g}} \Omega \cdot \nabla \rho_g+\mathcal{O}\left(\varepsilon^{2}\right).
\end{aligned}
\end{equation}
Taking the integral with respect to $\Omega$ of the equation \eqref{mgre1-equi-s0f},  one can get
\begin{equation} \label{diff-mgre-s0f}
\partial_t \rho_g  +\f{1}{\eps}\nabla \cdot \average{\Omega I_g} = \sigma_{a,g}(B_g-\rho_g)\,.
\end{equation}
Plugging \eqref{eqapp1-appen-s0f} into the equation \eqref{diff-mgre-s0f}, 
letting $\eps\to 0$, then \eqref{diff-mgre-s0f} reduces to
\begin{equation} \label{diff-mg1re-s0f}
\partial_t\rho_g-\nabla \cdot\left(  \frac{1}{3\sigma_{s,g}} \nabla \rho_g\right)=  \sigma_{a,g}(B_g-\rho_g)  \,.
\end{equation}
Moreover, by letting $\eps\to 0$ in \eqref{mgre2-equi-s0f}, one has
$$
 C_v\partial_t T  = 4\pi\sum\limits_{g'=1}^{G } \left[\sigma_{a,g'}\left(  \rho_{g'}-  B_{g'}\right)\right],
$$
and the above two equations give a group discretization of the frequency dependent diffusion limit \eqref{nonequ-diff}.
When one approximates $\sigma_{a,g}$ and $\sigma_{s,g}$ as in \eqref{eq:piecewisesection}, the coefficients in the frequency dependent diffusion limit are determined by
$$
\sigma^c_{a,g}=\f{1}{\nu_{g+1/2}-\nu_{g-1/2}}\int_{\nu_{g-\frac{1}{2}}}^{\nu_{g+\frac{1}{2}}} \sigma_{a}\mathrm{d} \nu,\quad\f{1}{\sigma_{s,g}^c}=\f{1}{\f{1}{\nu_{g+1/2}-\nu_{g-1/2}}\int_{\nu_{g-\frac{1}{2}}}^{\nu_{g+\frac{1}{2}}} \sigma_{s}\mathrm{d} \nu}
$$
When \eqref{eqn::ra} and \eqref{ros-s} are used, $\sigma^r_{a,g}$ and $1/\sigma_{s,g}^r$ become
$$
\sigma^r_{a,g}=\frac{\int_{\nu_{g-\frac{1}{2}}}^{\nu_{g+\frac{1}{2}}}\f{\sigma_a}{\sigma_s} \frac{\partial B}{\partial T}(T) \mathrm{~d} \nu}{\int_{\nu_{g-\frac{1}{2}}}^{\nu_{g+\frac{1}{2}}} \f{1}{\sigma_s} \frac{\partial B}{\partial T} (T)\mathrm{~d} \nu}, \quad\f{1}{\sigma^r_{s,g}}=\frac{\int_{\nu_{g-\frac{1}{2}}}^{\nu_{g+\frac{1}{2}}}\f{1}{\sigma_s} \frac{\partial B}{\partial T}(T) \mathrm{~d} \nu}{\int_{\nu_{g-\frac{1}{2}}}^{\nu_{g+\frac{1}{2}}} \frac{\partial B}{\partial T} (T)\mathrm{~d} \nu}.
$$
 We emphasize that, different approximations provide different coefficient values for different groups in the limiting models, which is crucial to get the correct solution behavior.
\subsection{Frequency discretization of the FDDL}\label{apfre} 
From the above discussion, we can see that different group discretizations for FRTE yield different approximations of the group diffusion and absorption coefficients in the FDDL. In the frequency dependent diffusion regime, different group diffusion and absorption coefficients can affect the propagation speed of heat waves under different physical conditions. Moreover, which approximation is better depends on different physical conditions \cite{yee2017stable,park2019multigroup}. 


	In this subsection,  we would like to compare several different approximations of the coefficients $\sigma$ in \eqref{eq:limit_T} and $\sigma_s$, $\sigma_a$ in \eqref{nonequ-diff}. We do not discuss about which coefficient approximation is better in the limiting model, but to show that they can be quite different from each other.

	Three different approximations of $\sigma_{a,g}$ are considered: the piece-wise constant approximation $\sigma^c_{a,g}$ in \eqref{eq:piecewisesection}, the Rosseland mean approximation $\sigma_{a,g}^r$ in \eqref{eqn::ra} and the Planck approximation. In the Planck approximation as in \cite{sunj}, the group absorptions to discretize \eqref{nonequ-diff} are defined as follows 
\begin{equation}
\sigma^e_{a, g}=\frac{\int_{\nu_{g-\frac{1}{2}}}^{\nu_{g+\frac{1}{2}}}\sigma_aB(\nu,T) \mathrm{~d} \nu}{\int_{\nu_{g-\frac{1}{2}}}^{\nu_{g+\frac{1}{2}}} B (\nu,T)\mathrm{~d} \nu}, \quad\sigma^a_{a, g}=\frac{\int_{\nu_{g-\frac{1}{2}}}^{\nu_{g+\frac{1}{2}}}\sigma_a\rho \mathrm{~d} \nu}{\int_{\nu_{g-\frac{1}{2}}}^{\nu_{g+\frac{1}{2}}} \rho	\mathrm{~d} \nu}.
\end{equation}
It should be pointed out that $\sigma^a_{a, g}$ is a weighted integration with the unknown function $\rho$. Usually,
the unknown function $\rho$ is replaced by the Planck function with the radiation temperature $T_r$, and the multi-group FDDL writes:
\begin{subequations}\label{eqn:003re-tr-p}  
	\begin{numcases}{}
	\partial_t \rho_g+\nabla\cdot\left(\frac{1}{3\sigma_{s,g}} \nabla \rho_g\right) =\sigma^e_{a,g} B_g-\sigma^a_{a,g}\rho_g,\qquad (g=1,\cdots,G),  \label{mgre1-tr-p}\\
	C_v\partial_t T  =  4\pi\sum\limits_{g'=1}^{G } \left( \sigma^a_{a,g'} \rho_{g'}- \sigma^e_{a,g'} B_{g'}\right).\label{mgre2-tr-p}
	\end{numcases}
\end{subequations}

	Similar as in \cite{densmore2012hybrid,sunj}, we consider a simple dependence of the frequency such that $\sigma_a=\nu^{-3}$ and $\sigma_s=\nu^{-3}$. Here we ignore the scale parameters $\mathscr{L}_a$ and $\mathscr{L}_s$ since they not affect the results of the comparison. The frequency domain is from $0.1\ eV$ to $100\ keV$ and divided into 30 frequency groups, and the groups are logarithmically spaced.

\begin{figure}[!h]
	\centering
		\subfigure[]{
		\includegraphics[width=3.2in]{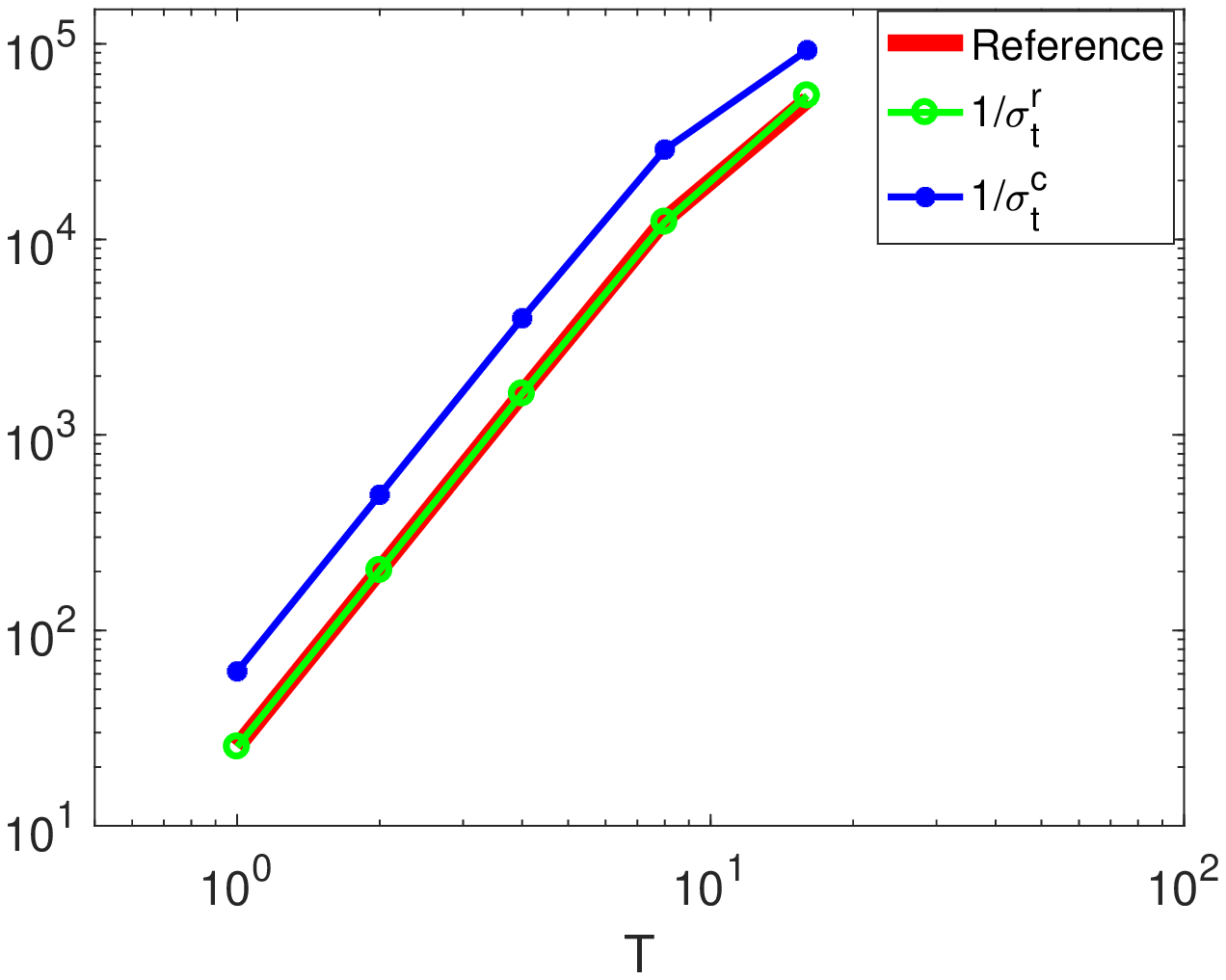}
		\label{fig:striped_diff}
			}
		\subfigure[]{
		\includegraphics[width=3.2in]{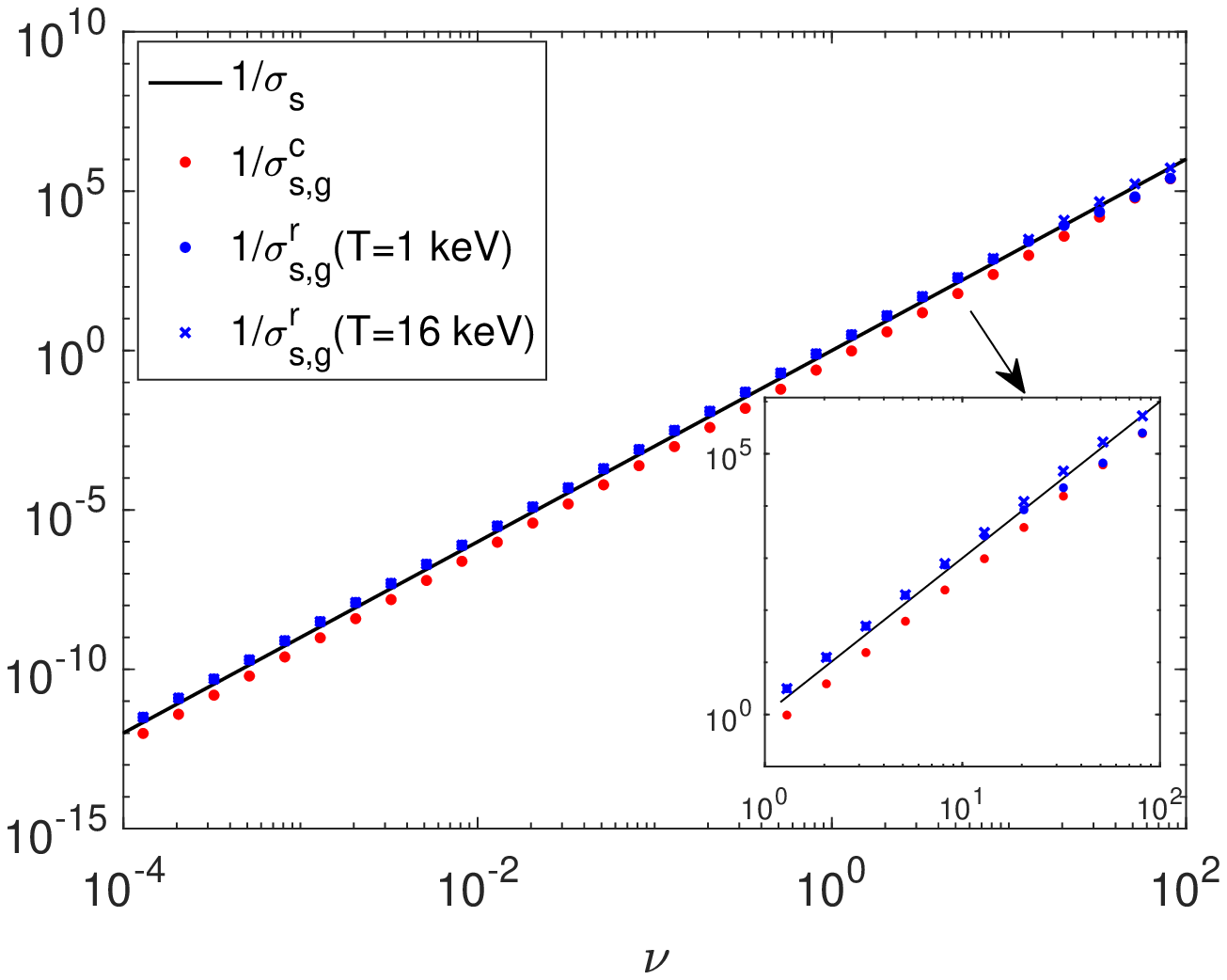}
		\label{fig.2_ssg}
	}
\caption{Left: comparison of the different approximations of the mean free path  $1/\sigma$ in the gray radiation diffusion limit \eqref{eq:defsigma} at the different temperatures; Right: comparison of the different approximations of the diffusion coefficient $1/\sigma_s$ in the FDDL \eqref{nonequ-diff}    at  $T=1\ keV$ and $T=16\ keV$.}
	\label{fig:striped_diff-1}
\end{figure}

\begin{figure}[!h]
	\centering
		\subfigure[]{
		\includegraphics[width=3.22in]{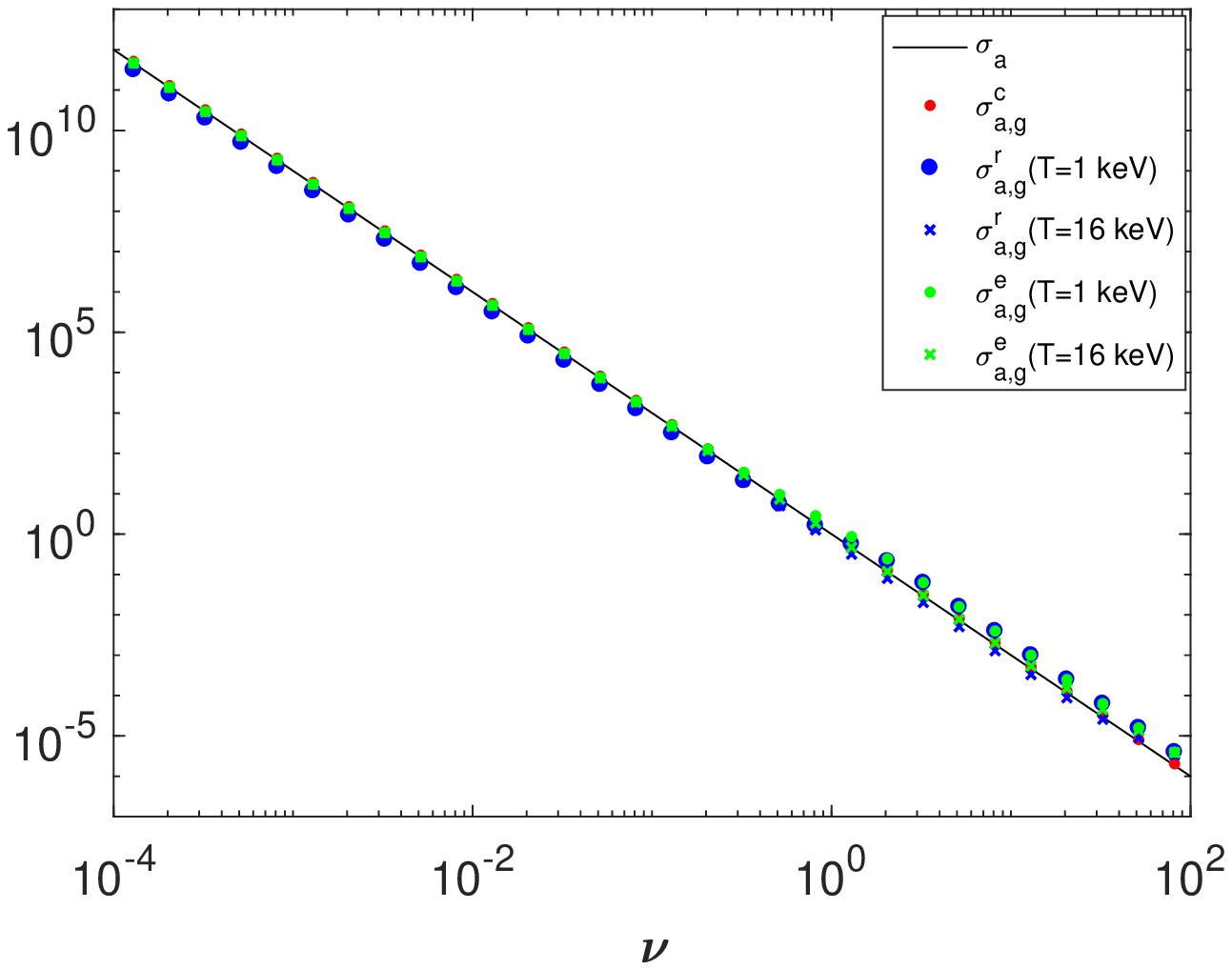}
	}
	\subfigure[]{
		\includegraphics[width=3.2in]{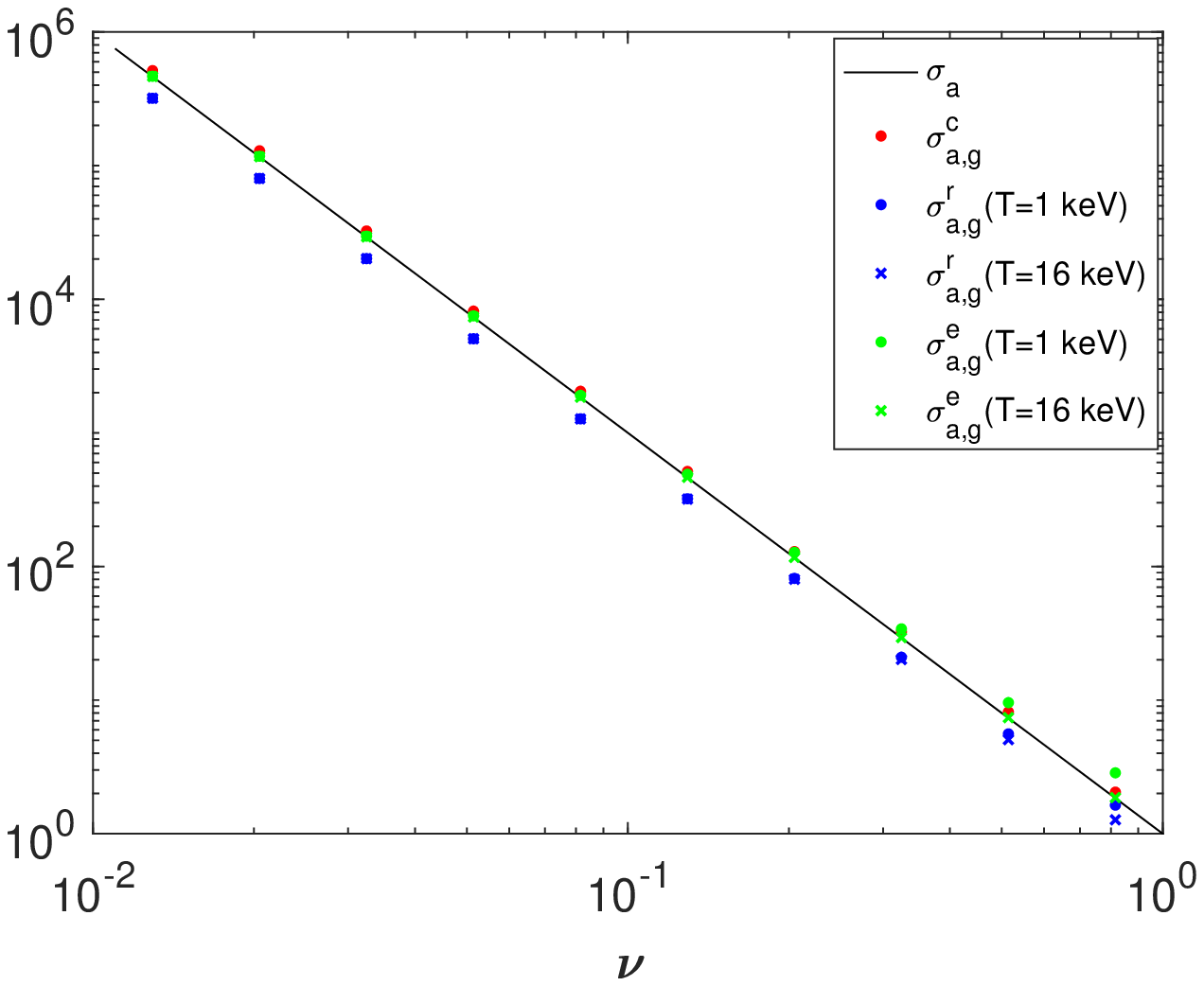}
	}\\
		\subfigure[]{
		\includegraphics[width=3.2in]{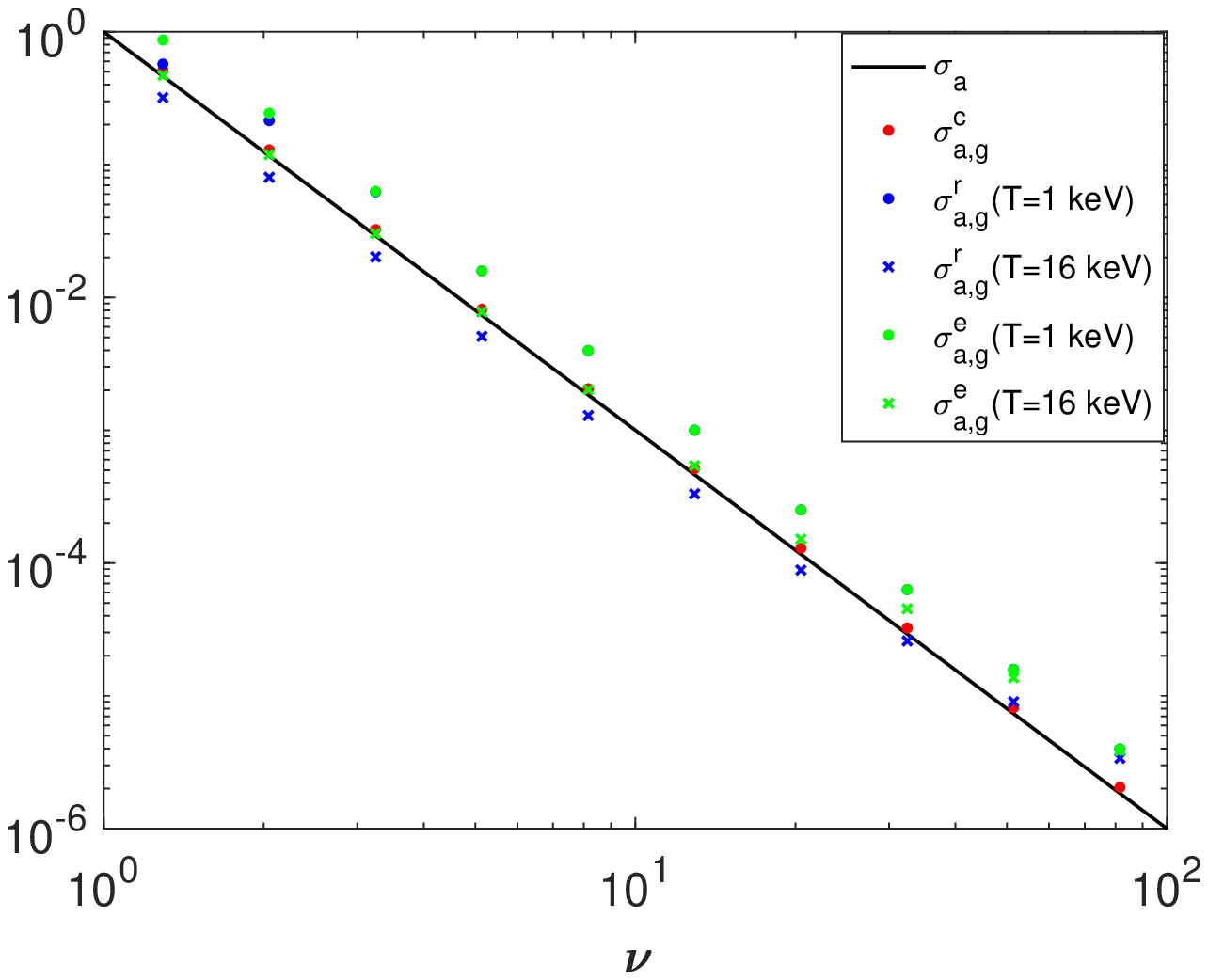}
	}\\
	\caption{(a)comparison of the different approximations of  $\sigma_a$  at $T= 1\ keV$ and $T= 16\ keV$; (b)zoom in part of figure (a), here the some red dots and green dots overlap; (c)zoom in, here the blue dots and green dots overlap. }
	\label{figequ-2sa}
\end{figure}

The following quantities are compared:
\begin{itemize}
\item The two different approximations of $1/\sigma_t$ in \eqref{sigmatconstant} and \eqref{eq:sigmat1} at different temperatures $T$. The results are displayed in Table~\ref{tab:stability1} and Figure.~\ref{fig:striped_diff}, and the reference solution is calculated by dividing the frequency domain into 600 groups.
\item  The diffusion coefficients $1/\sigma_{s}$ in the FDDL \eqref{nonequ-diff} are approximated by the $1/\sigma^c_{s,g}$ and two different $1/\sigma^r_{s,g}$ at $T=1\ keV$ and $T=16\ keV$. In Figure.~\ref{fig.2_ssg}, From Figure.~\ref{fig:striped_diff-1}, the Rosseland mean approximation is closer to the reference. 
\item For the absorbing coefficient $\sigma_a$, we compare the piece-wise constant approximation $\sigma^c_{a,g}$, the Rosseland mean approximation $\sigma_{a,g}^r$ and the Plank approximation $\sigma^e_{a,g}$ in Figure.~\ref{figequ-2sa}. We observe that they can be very different.
\end{itemize}
\begin{table}[!ht]
	\centering
	\caption{The different approximations to  $1/\sigma$ in the gray radiation diffusion limit \eqref{eq:defsigma}. The reference solution is calculated by dividing the frequency domain into 600 groups.}
	\vspace{0.3cm}
	\begin{tabular}{l|ccccc}
		\hline
		T&Reference &$1/\sigma_t^r$ &relative error& $1/\sigma_t^c$&relative error \\
		\hline
		$1\ keV$ &25.3901&25.3920&$7.48\times10^{-5}$&	61.7458& 1.43\\
		$2\ keV$ &203.1204&203.1042&$7.97\times10^{-5}$&	493.8879&1.43  \\
		$4\ keV$ &	$1.6250\times10^3$&$1.6251\times10^3$&$6.15\times10^{-5}$&	$3.9516\times10^3$&1.43 \\ 
		$8\ keV$ &$1.2574\times10^4$&$1.2284\times10^4$&$2.31\times10^{-2}$&	$2.8811\times10^4$&1.29  \\ 
		$16\ keV$ &$5.1674\times10^4$&$5.4364\times10^4$&$5.21\times10^{-2}$&	$9.2914\times10^4$&$7.98\times10^{-1}$  \\ 
		\hline
	\end{tabular}\vspace{0cm}
	\label{tab:stability1}
\end{table}



\section{The decomposed multi-group method and its asymptotic limits}
\subsection{The decomposed multi-group method}
A new decomposed multi-group method is proposed in this section. We decompose the radiation intensity $I(t,x,\Omega,\nu)$ into three parts:
 \begin{subequations}\label{eqn:011}
 \begin{equation}\label{eqn:0111}
 I(t,x,\Omega,\nu)=\average{I}+3 \Omega\cdot\average{\Omega I}+ Q(t,x,\Omega,\nu),
 \end{equation}
  \begin{equation} \label{eqn:0112}
\hspace{3.2cm}:=\rho(t,x,\nu)+ \Omega\cdot R(t,x,\nu)+ Q(t,x,\Omega,\nu),
 \end{equation}
 \end{subequations}
with $\average{\cdot}=\f{1}{4\pi}\int_{4\pi} \cdot \ \text{d}\Omega$ and 3$\average{\Omega I}= R$. Then taking $\average{\Omega\cdot}$  on both sides of \eqref{eqn:0111} and using the condition $\average{\Omega\Omega}:=D_d=\f{1}{3}I_d$  (where $I_d$ denotes the 3 by 3 identity matrix), one gets
\begin{equation*}
\average{\Omega I}=\average{\Omega I}+ \average{\Omega Q},
\end{equation*}
which indicates that
\begin{equation}\label{eqn:rr}
\average{\Omega Q}=0.
\end{equation}
Moreover,  taking $\average{\cdot}$  on both sides of \eqref{eqn:011}, one gets
$$
\average{I}=\average{\rho+ \Omega\cdot R+ Q}=\average{I}+\average{ Q},
$$
which implies
\begin{equation}\label{eqn:qq}
\average{ Q}=0.
\end{equation}
\eqref{eqn:rr} and \eqref{eqn:qq} are two properties that $Q$ has to satisfy.

We reconstruct the intensity $I$ in each interval $(\nu_{g-1/2},\nu_{g+1/2})$ in the frequency domain as follows:
\begin{equation}\label{linearreconstruction}
I(t,x,\Omega,\nu)=\frac{1}{\nu_{g+1/2}-\nu_{g-1/2}}\textcolor{black}{\Big(\omega_1\rho_g+\Omega\cdot R_g+Q_g+s_{\rho_g}+s_{\Omega\cdot R_g} \Big)}
\end{equation}
where
$$
\rho_g=\average{I_g},\qquad
R_g=3\average{\Omega I_g},\qquad Q_g=I_g-\rho_g-\Omega\cdot R_g,
$$
and
\textcolor{black}{
\begin{equation}\label{eqn:slope}
	s_{\rho_g}=\textcolor{black}{(\nu_{g+1/2}-\nu_{g-1/2})B(\nu,T)-\omega_2B_g},\qquad s_{\Omega\cdot R_g}=\f{\f{\nu_{g+1/2}-\nu_{g-1/2}}{\mathscr{L}_a\sigma_a+\mathscr{L}_s\sigma_s}\f{\partial B}{\partial T}}{\int_{\nu_{g-\frac{1}{2}}}^{\nu_{g+\frac{1}{2}}}\f{1}{\mathscr{L}_a\sigma_a+\mathscr{L}_s\sigma_s}\f{\partial B}{\partial T} \text{d} \nu
		}\Omega\cdot R_g-\Omega\cdot R_g.
\end{equation}}
Here $\omega_1$ and $\omega_2$ are two given weight functions that depend on $\nu$. We will see later on that by tuning $\omega_1$, $\omega_2$, one can get the required approximations of the group absorption coefficients in FDDL. 
The first three terms on the right hand side of \eqref{linearreconstruction} $\rho_g$, $\Omega\cdot R_g$, $Q_g$ depend only on $g$ not in $\nu$. Moreover, the zeroth and first moments of $Q_g$ equal to 0, which is the same as $Q$. The two additional terms $s_{\rho_g}$, $s_{\Omega\cdot R_g}$ depend on $\nu$ but satisfy $\int\limits_{\nu_{g-1/2}}^{\nu_{g+1/2}}s_{\rho_g}\text{d} \nu=\int\limits_{\nu_{g-1/2}}^{\nu_{g+1/2}}s_{\Omega\cdot R_g}\text{d} \nu=0$. Thus the reconstructed $I$ in \eqref{linearreconstruction} satisfies $\int\limits_{\nu_{g-1/2}}^{\nu_{g+1/2}}I\text{d} \nu=I_g$. 
 The two additional terms $s_{\rho_g}$ and $s_{\Omega\cdot R_g}$ are important to get the correct limits in the gray radiation and frequency dependent diffusion regime. The detailed derivations will be displayed in section 3.3. 

It is important to note that, as far as $I_g, T$ are given, $\rho_g$, $R_g$, $Q_g$ and $s_{\rho_g}$, $s_{\Omega\cdot R_g}$ are known.
 Then $I(t,x,\Omega,\nu)$ in \eqref{linearreconstruction} can be determined by $I_g$ and the two terms in \eqref{eq:sigmaasiamgsdecomposed} are approximated by
$$
\begin{aligned}
&\int\limits_{\nu_{g-1/2}}^{\nu_{g+1/2}}\sigma_a\left( B-I\right)\text{d}\nu=\int\limits_{\nu_{g-\frac{1}{2}}}^{\nu_{g+\frac{1}{2}}}\sigma_a\Big[B-\frac{1}{\nu_{g+1/2}-\nu_{g-1/2}}\textcolor{black}{\big(\omega_1\rho_g+\Omega\cdot R_g+Q_g+s_{\rho_g}+s_{\Omega\cdot R_g} \big)} \Big] \text{d} \nu\\
&\hspace{-0.6cm}=\int\limits_{\nu_{g-\frac{1}{2}}}^{\nu_{g+\frac{1}{2}}}\sigma_a\Big[B-\frac{1}{\nu_{g+1/2}-\nu_{g-1/2}}\textcolor{black}{\Big(\omega_1\rho_g+(\nu_{g+1/2}-\nu_{g-1/2})B-\omega_2B_g+\f{\f{\nu_{g+1/2}-\nu_{g-1/2}}{\mathscr{L}_a\sigma_a+\mathscr{L}_s\sigma_s}\f{\partial B}{\partial T}}{\int_{\nu_{g-\frac{1}{2}}}^{\nu_{g+\frac{1}{2}}}\f{1}{\mathscr{L}_a\sigma_a+\mathscr{L}_s\sigma_s}\f{\partial B}{\partial T} \text{d} \nu
}\Omega\cdot R_g+ Q_g\Big)}\Big] \text{d} \nu\\
&\hspace{-0.6cm}\textcolor{black}{=\frac{1}{\nu_{g+1/2}-\nu_{g-1/2}}\int\limits_{\nu_{g-\frac{1}{2}}}^{\nu_{g+\frac{1}{2}}}\sigma_a(\omega_2B_g-\omega_1\rho_g)\text{d} \nu-\frac{1}{\nu_{g+1/2}-\nu_{g-1/2}}\int\limits_{\nu_{g-\frac{1}{2}}}^{\nu_{g+\frac{1}{2}}}\sigma_a\text{d} \nu Q_g}\\
&\hspace{7cm}-\int\limits_{\nu_{g-\frac{1}{2}}}^{\nu_{g+\frac{1}{2}}}\f{\sigma_a}{\mathscr{L}_a\sigma_a+\mathscr{L}_s\sigma_s}\f{\partial B}{\partial T}\text{d} \nu\f{\Omega\cdot R_g}{\int_{\nu_{g-\frac{1}{2}}}^{\nu_{g+\frac{1}{2}}}\f{1}{\mathscr{L}_a\sigma_a+\mathscr{L}_s\sigma_s}\f{\partial B}{\partial T} \text{d} \nu
},
\end{aligned}
$$
and
$$
\begin{aligned}
\int\limits_{\nu_{g-1/2}}^{\nu_{g+1/2}}\sigma_s\left(\rho-I\right)\text{d}\nu&=\frac{-1}{\nu_{g+1/2}-\nu_{g-1/2}}\int\limits_{\nu_{g-\frac{1}{2}}}^{\nu_{g+\frac{1}{2}}}\sigma_s\big(\Omega \cdot R_g +Q_g+s_{\Omega\cdot R_g}\big)\text{d} \nu\\
&=\frac{-1}{\nu_{g+1/2}-\nu_{g-1/2}}\int\limits_{\nu_{g-\frac{1}{2}}}^{\nu_{g+\frac{1}{2}}}\sigma_s\text{d} \nu Q_g-\int\limits_{\nu_{g-\frac{1}{2}}}^{\nu_{g+\frac{1}{2}}}\f{\sigma_s}{\mathscr{L}_a\sigma_a+\mathscr{L}_s\sigma_s}\f{\partial B}{\partial T}\text{d} \nu\f{\Omega\cdot R_g}{\int_{\nu_{g-\frac{1}{2}}}^{\nu_{g+\frac{1}{2}}}\f{1}{\mathscr{L}_a\sigma_a+\mathscr{L}_s\sigma_s}\f{\partial B}{\partial T} \text{d} \nu
}.
\end{aligned}
$$
Then we can get the decomposed multi-group radiative transfer equation:
	\begin{subequations}\label{eqn:003re-equi-z}  
		\begin{numcases}{}
			\frac{1}{\mathcal{C}}\partial_t I_g+\Omega \cdot \nabla I_g =\mathscr{L}_a \sigma^2_{a,g} B_g- \mathscr{L}_a \sigma^1_{a,g}\rho_g-\big(\mathscr{L}_a+\mathscr{L}_s\big)\sigma_{t,g}\Omega\cdot R_g-(\mathscr{L}_a\sigma_{a,g}+\mathscr{L}_s\sigma_{s,g}) Q_g,\nonumber\\
			\hspace{12cm}(g=1,\cdots,G),  \label{mgre1-equi-z}\\
			 C_v\partial_t T  = 4\pi \mathcal{C}\PP_0\sum\limits_{g'=1}^{G } \left[\mathscr{L}_a \left(  \sigma^1_{a,g'} \rho_{g'}-   \sigma^2_{a,g'}B_{g'}\right)\right],\label{mgre2-equi-z}
		\end{numcases}
	\end{subequations}
	where
	\begin{equation}\label{eq:groupcoefficient}
	\begin{aligned}
		&\sigma^1_{a,g}=\frac{1}{\nu_{g+1/2}-\nu_{g-1/2}}\int\limits_{\nu_{g-1/2}}^{\nu_{g+1/2}}\omega_1\sigma_{a}\text{d}\nu,\quad\sigma^2_{a,g}=\frac{1}{\nu_{g+1/2}-\nu_{g-1/2}}\int\limits_{\nu_{g-1/2}}^{\nu_{g+1/2}}\omega_2\sigma_{a}\text{d}\nu,\\
	&\sigma_{a,g}=\frac{1}{\nu_{g+1/2}-\nu_{g-1/2}}\int\limits_{\nu_{g-1/2}}^{\nu_{g+1/2}}\sigma_{a}\text{d}\nu,\quad\sigma_{s,g}=\frac{1}{\nu_{g+1/2}-\nu_{g-1/2}}\int\limits_{\nu_{g-1/2}}^{\nu_{g+1/2}}\sigma_{s}\text{d}\nu,\\
	&\hspace{3cm}\sigma_{t,g}=\f{1/(\mathscr{L}_a+\mathscr{L}_s)}{\int_{\nu_{g-\frac{1}{2}}}^{\nu_{g+\frac{1}{2}}}\f{1}{\mathscr{L}_a\sigma_{a}+\mathscr{L}_s \sigma_s} \f{\partial B}{\partial T}\text{d} \nu\big/\int_{\nu_{g-\frac{1}{2}}}^{\nu_{g+\frac{1}{2}}}\f{\partial B}{\partial T} \text{d} \nu}.
	\end{aligned}
	\end{equation}
 
 It is important to use the reconstruction as in \eqref{linearreconstruction} to get the required absorption coefficients for group discretizations of the FDDL. One can consider the following possible $\omega_1$ and $\omega_2$: 
\begin{itemize}
	\item When $\omega_1=\omega_2=1$, then the absorbing coefficients are
	$$
	\sigma^1_{a,g}=\sigma^2_{a,g}=\frac{1}{\nu_{g+1/2}-\nu_{g-1/2}}\int\limits_{\nu_{g-1/2}}^{\nu_{g+1/2}}\sigma_{a}\text{d}\nu,
	$$
	which will yield piece-constant approximations of the  group absorption coefficients in FDDL.
	\item When $\omega_1=\f{\f{\partial B}{\partial T}(T_r)}{\frac{1}{\nu_{g+1/2}-\nu_{g-1/2}}\f{\partial B_g}{\partial T}(T_r)},\quad\omega_2=\f{\f{\partial B}{\partial T}(T)}{\frac{1}{\nu_{g+1/2}-\nu_{g-1/2}}\f{\partial B_g}{\partial T}(T)}$, then the absorbing coefficients are
	$$
	\sigma^1_{a,g}=\f{\int_{\nu_{g-1/2}}^{\nu_{g+1/2}}\sigma_{a}\f{\partial B}{\partial T}(T_r)\text{d}\nu}{\f{\partial B_g}{\partial T}(T_r)},\quad
	\sigma^2_{a,g}=\f{\int_{\nu_{g-1/2}}^{\nu_{g+1/2}}\sigma_{a}\f{\partial B}{\partial T}(T)\text{d}\nu}{\f{\partial B_g}{\partial T}(T)},
	$$
	which gives the Rosseland mean approximations of the group absorption coefficients in FDDL.
	\item When  $\omega_1=\f{ B(T_r)}{\frac{1}{\nu_{g+1/2}-\nu_{g-1/2}} B_g(T_r)},\quad\omega_2=\f{ B(T)}{\frac{1}{\nu_{g+1/2}-\nu_{g-1/2}} B_g(T)}$, then the absorbing coefficients are
	$$
	\sigma^1_{a,g}=\f{\int_{\nu_{g-1/2}}^{\nu_{g+1/2}}\sigma_{a}B(T_r)\text{d}\nu}{ B_g(T_r)},\quad
	\sigma^2_{a,g}=\f{\int_{\nu_{g-1/2}}^{\nu_{g+1/2}}\sigma_{a}B(T)\text{d}\nu}{ B_g(T)},
	$$
	which corresponds to the Planck approximation of the  absorption coefficient in FDDL.
\end{itemize}	
When the different $\omega_1$ and $\omega_2$ are chosen, the different approximations of the absorbing coefficients in FDDL are determined.
Comparing  \eqref{eqn:003re-equi-z} with \eqref{eqn:003re-tr}, we can see that the scattering and absorbing coefficients for $\rho_g$, $R_g$ and $Q_g$ are different in the decomposed multi-group method, while $\rho_g$, $R_g$ and $Q_g$ share the same values of scattering and absorbing coefficients in the classical multi-group method. 

\subsection{The decomposed multi-group method in the diffusion regime}
\subsubsection{The gray radiation diffusion regime}
For the scaling $\mathscr{L}_a=1/\eps,  \mathscr{L}_s=1/\eps$, the decomposed  multi-group radiative transfer system \eqref{eqn:003re-equi-z}   can be rewritten into \textcolor{black}{
\begin{subequations}\label{eqn:003re-equi}  
\begin{numcases}{}
\eps\partial_t I_g+\Omega \cdot \nabla I_g = \f{
1}{\eps} \left( \sigma_{a,g}^2 B_g-\sigma_{a,g}^1 \rho_g\right)-\f{\Omega\cdot R_g}{\int_{\nu_{g-\frac{1}{2}}}^{\nu_{g+\frac{1}{2}}}\f{\eps}{ \sigma_{a}+ \sigma_s} \f{\partial B}{\partial T}\text{d} \nu\big/\int_{\nu_{g-\frac{1}{2}}}^{\nu_{g+\frac{1}{2}}} \f{\partial B}{\partial T}\text{d} \nu
}-\f{ \sigma_{a,g}+ \sigma_{s,g}}{\eps} Q_g,\nonumber\\
 \hspace{10.5cm}(g=1,\cdots,G),  \label{mgre1-equi}\\
\eps C_v\partial_t T  = 4\pi\sum\limits_{g'=1}^{G } \left[\f{1}{\eps}\left( \sigma^1_{a,g'} \rho_{g'}-  \sigma^2_{a,g'}B_{g'}\right)\right].\label{mgre2-equi}
\end{numcases}
\end{subequations}
At the leading order of \eqref{mgre1-equi}, one has
\begin{align}
\sigma^2_{a,g}B_g-\sigma^1_{a,g}\rho_g- \f{\Omega\cdot R_g}{\int_{\nu_{g-\frac{1}{2}}}^{\nu_{g+\frac{1}{2}}}\f{1}{\sigma_{a}+\sigma_{s} }\f{\partial B}{\partial T} \text{d} \nu\big/\int_{\nu_{g-\frac{1}{2}}}^{\nu_{g+\frac{1}{2}}} \f{\partial B}{\partial T}\text{d} \nu
}-(\sigma_{a,g}+\sigma_{s,g} )Q_g=\order(\eps),\label{anatz1re}
\end{align}
Integrating with respect to $\Omega$ in equation \eqref{anatz1re} and using $\average{Q_g}=0$, we can obtain
\begin{equation}\label{rhogBg-zj}
\sigma^1_{a,g}\rho_g = \sigma^2_{a,g} B_g + \order(\eps)\,.
\end{equation}
Summing the above equation over the group index $g$, and by the definitions \eqref{eqapp}, one can get
\begin{equation}
\int\limits_{0}^\infty\sigma_a\rho\text{d}\nu=\int\limits_{0}^\infty\sigma_aB\text{d}\nu + \order(\eps),
\end{equation}
which implies that
\begin{equation}\label{rhogBg}
\rho_g = B_g + \order(\eps)\,.
\end{equation}}
Integrating \eqref{mgre1-equi} with respect to $\Omega$ and adding up all groups and equation \eqref{mgre2-equi}, one gets
\begin{equation} \label{diff-mgre}
	4\pi\partial_t\left(\sum_{g=1}^G \rho_g\right)  +\f{4\pi}{3\eps} \sum_{g=1}^G\nabla \cdot R_g + C_v \partial_t T = 0\,.
\end{equation}
From \eqref{rhogBg} and the conditions that $B_g$ satisfies in \eqref{Bgconditions}, $4\pi\sum_{g=1}^G \rho_g= T^4+\order(\eps)$. In order 
to get the diffusion limit equation, one has to determine the relationship between $\rho_g$ and $R_g$. Taking the first moment of \eqref{mgre1-equi}, and using $I_g=\rho_g+\Omega\cdot R_g+Q_g$, \eqref{rhogBg} and $\average{\Omega Q_g}$=0, one has
\begin{equation}\label{eqnrg---}
	-\f{1}{\eps}R_g=\frac{1}{\int_{\nu_{g-\frac{1}{2}}}^{\nu_{g+\frac{1}{2}}} \f{\partial B}{\partial T}\text{d} \nu}\left[\int_{\nu_{g-\frac{1}{2}}}^{\nu_{g+\frac{1}{2}}}\f{\nabla\rho_g}{\sigma_{a}+\sigma_{s} }\f{\partial B}{\partial T} \text{d} \nu+\int_{\nu_{g-\frac{1}{2}}}^{\nu_{g+\frac{1}{2}}}\f{\nabla\average{\Omega\Omega Q_g}}{\sigma_{a}+\sigma_s } \f{\partial B}{\partial T}\text{d} \nu\right] +\order(\eps).
\end{equation}
Taking the second moment of $\Omega$ in equation \eqref{anatz1re} and noting \eqref{rhogBg} yields
$$
	\average{\Omega\Omega Q_g}:=K_{Q_g}=\order(\eps),
$$
which implies
\begin{equation}\label{eqnrg}
	-\f{1}{\eps}R_g=\frac{1}{\int_{\nu_{g-\frac{1}{2}}}^{\nu_{g+\frac{1}{2}}} \f{\partial B}{\partial T}\text{d} \nu}\int_{\nu_{g-\frac{1}{2}}}^{\nu_{g+\frac{1}{2}}}\f{\nabla\rho_g}{\sigma_{a}+\sigma_{s}  }\f{\partial B}{\partial T} \text{d} \nu+\order(\eps).
\end{equation}
Plugging \eqref{rhogBg}, \eqref{eqnrg} into the equation \eqref{diff-mgre},
when $\eps\to 0$, \eqref{diff-mgre} reduces to
\begin{equation} \label{diff-mg1re}
 \partial_t\left(\sum_{g=1}^G 4\pi B_g\right) + C_v\partial_t T=  \sum_{g=1}^G\nabla \cdot\left[ \left(\int\limits_{\nu_{g-1/2}}^{\nu_{g+1/2}}  \frac{4\pi}{3(\sigma_a+\sigma_s)} \f{\partial B}{\partial T}\text{d} \nu\right)\nabla T\right]  \,.
\end{equation}
Due to \eqref{Bgconditions}, \eqref{diff-mg1re} becomes
\begin{equation} \label{diff-mg2}
\partial_t T^4  + C_v\partial_t T=  \nabla \cdot \left(  \frac{1}{3\widehat\sigma_R} \nabla T^4 \right) \,,
\end{equation}
where
$$
\f{1}{\widehat\sigma_R}=\f{1}{4T^3} \left(\sum\limits_{g=1}^{G }\int\limits_{\nu_{g-1/2}}^{\nu_{g+1/2}}  \frac{4\pi}{\sigma_a+\sigma_s} \f{\partial B}{\partial T}\text{d} \nu\right),
$$
which is consistent with equation \eqref{eq:limit_T}. 
\subsubsection{The frequency dependent diffusion regime}\label{nedl}
In the frequency dependent diffusion regime, the scaling is $\mathscr{L}_a=\eps$  and $\mathscr{L}_s=1/\eps$, then the multi-group radiative transfer system \eqref{eqn:003re-equi-z} can be rewritten to \textcolor{black}{
\begin{subequations}\label{eqn:003re-nonequi}  
\begin{numcases}{}
	\eps\partial_t I_g+\Omega \cdot \nabla I_g = \eps\left(  \sigma^2_{a,g}B_g- \sigma^1_{a,g}\rho_g\right)-\f{\Omega\cdot R_g}{\int_{\nu_{g-\frac{1}{2}}}^{\nu_{g+\frac{1}{2}}}\f{1}{\eps \sigma_a+\f{\sigma_{s}}{\eps}}\f{\partial B}{\partial T} \text{d} \nu\big/\int_{\nu_{g-\frac{1}{2}}}^{\nu_{g+\frac{1}{2}}} \f{\partial B}{\partial T}\text{d} \nu
	}-\left(\eps \sigma_{a,g}+\f{\sigma_{s,g}}{\eps}\right) Q_g,\nonumber\\
	\hspace{12cm}(g=1,\cdots,G),   \label{mgre1-nonequi}\\
	\eps C_v\partial_t T  = 4\pi\sum\limits_{g'=1}^{G } \left[\eps \left(\sigma^1_{a,g'}  \rho_{g'}- \sigma^2_{a,g'} B_{g'}\right)\right].\label{mgre2-nonequi}
\end{numcases}
\end{subequations}}
The asymptotic analysis for the FDDL begins to find the relation between $\rho_g$ and $R_g$. From \eqref{mgre1-nonequi}, one has
\begin{align}
- \f{\Omega\cdot R_g}{\int_{\nu_{g-\frac{1}{2}}}^{\nu_{g+\frac{1}{2}}}\f{1}{\sigma_{s}}\f{\partial B}{\partial T} \text{d} \nu\big/\int_{\nu_{g-\frac{1}{2}}}^{\nu_{g+\frac{1}{2}}} \f{\partial B}{\partial T}\text{d} \nu
	}-\sigma_{s,g} Q_g=\order(\eps).\label{anatz1re-n}
\end{align}
Taking the second moment in $\Omega$ of equation \eqref{anatz1re-n} yields
\begin{equation}\label{kqg-n}
K_{Q_g}=\order(\eps).
\end{equation}
Moreover, using the relation \eqref{kqg-n}, the leading order of the first moment of  equation \eqref{mgre1-nonequi} is
\begin{equation}\label{eqnrg-n}
	-\f{1}{\eps}R_g=\frac{1}{\int_{\nu_{g-\frac{1}{2}}}^{\nu_{g+\frac{1}{2}}} \f{\partial B}{\partial T}\text{d} \nu}\int_{\nu_{g-\frac{1}{2}}}^{\nu_{g+\frac{1}{2}}}\f{\nabla\rho_g}{\sigma_{s} }\f{\partial B}{\partial T} \text{d} \nu+\order(\eps),
\end{equation}
which gives the relation between $R_g$ and $\rho_g$. Integrating \eqref{mgre1-nonequi} with respect to $\Omega$, one can get
\begin{equation} \label{diff-nonmgre}
 \partial_t \rho_g + \frac{1}{3\eps}  \nabla\cdot R_g  = \sigma^2_{a,g}B_g-\sigma^1_{a,g}\rho_g\,.
\end{equation}
Plugging \eqref{eqnrg-n} into \eqref{diff-nonmgre}, and letting $\eps\to 0$, \eqref{diff-nonmgre} gives 
\begin{equation} \label{diff-nonmgre1}
\partial_t\rho_g  - \nabla \cdot \left(\f{1}{3\sigma'_{s,g} }\nabla   \rho_g \right)  =\sigma^2_{a,g}B_g-\sigma^1_{a,g}\rho_g\,,
\end{equation}
where $\f{1}{\sigma'_{s,g}}=\int_{\nu_{g-\frac{1}{2}}}^{\nu_{g+\frac{1}{2}}}\f{1}{\sigma_{s} }\f{\partial B}{\partial T} \text{d} \nu\big/\int_{\nu_{g-\frac{1}{2}}}^{\nu_{g+\frac{1}{2}}}\f{\partial B}{\partial T} \text{d} \nu.$
From equation \eqref{mgre2-nonequi}, one has
$$
 C_v\partial_t T  =  \sum\limits_{g'=1}^{G } \left[4\pi(\sigma^1_{a,g'}\rho_{g'}-\sigma^2_{a,g'}B_{g'})\right].
$$
In summary, the diffusion limit of the decomposed multi-group radiative transfer system \eqref{eqn:003re-nonequi}
is
\begin{subequations}\label{eqn:003re-nonequi-diff}  
\begin{numcases}{}
\partial_t\rho_g  - \nabla \cdot \left(\f{1}{3\sigma'_{s,g} }\nabla   \rho_g \right)  =\sigma^2_{a,g}B_g-\sigma^1_{a,g}\rho_g\,,  \label{mgre1-nonequi-diff}\\
 C_v\partial_t T  =  \sum\limits_{g'=1}^{G } \left[4\pi(\sigma^1_{a,g'}\rho_{g'}-\sigma^2_{a,g'}B_{g'})\right].\label{mgre2-nonequi-diff}
\end{numcases}
\end{subequations}
The decomposed multi-group radiative transfer equations \eqref{eqn:003re} is AP in both gray radiation diffusion limit and the FDDL. Moreover, one can choose different approximations of $\sigma_{a,g}^1$, $\sigma_{a,g}^2$ in the limiting frequency discretization.
\section{AP time discretization of the decomposed multi-group method}
\subsection{The time discretization}
Based on the decomposition in \eqref{eqn:011} for radiation intensity $I$ and the decomposed multi-group method in \eqref{eqn:003re-equi-z}, we will present a time discretization for the FRTE and show its AP property.

Let
$$
t^n=n\Delta t, \qquad n=0,1,2,\cdots,
$$
and
$$
I^n\approx I(t^n,x,\Omega,\nu),\qquad \rho^n\approx \rho(t^n,x,\nu),\qquad T^n\approx T(t^n,x),\qquad B^n\approx B(\nu,T(t^n,x)),
$$
with
$$
I_g^n\approx \int\limits_{\nu_{g-\frac{1}{2}}}^{\nu_{g+\frac{1}{2}}}I(t^n,x,\Omega,\nu) \text{d} \nu,\qquad \rho_g^n\approx \int\limits_{\nu_{g-\frac{1}{2}}}^{\nu_{g+\frac{1}{2}}} \rho(t^n,x,\nu)\text{d} \nu,\qquad B_g^n\approx \int\limits_{\nu_{g-\frac{1}{2}}}^{\nu_{g+\frac{1}{2}}} B(\nu,T(t^n,x))\text{d} \nu.
$$
We use the following semi-discrete scheme in time for \eqref{eqn:003re} in each group:
\begin{subequations}\label{eqn:003re-dis}  
	\begin{numcases}{}
		\f{1}{\mathcal{C}} \f{I_g^{n+1}-I_g^{n}}{\Delta t}+\Omega \cdot \nabla (\rho^{n+1}_g+\Omega\cdot R^{n+1}_g+Q^n_g )  =\nonumber\\ \hspace{0.5cm}\mathscr{L}_a \left( \sigma^2_{a,g}B^{n+1}_g-\sigma^1_{a,g} \rho^{n+1}_g\right)-(\mathscr{L}_a+\mathscr{L}_s)\sigma_{t,g}\Omega\cdot R^{n+1}_g-(\mathscr{L}_a\sigma_{a,g}+\mathscr{L}_s\sigma_{s,g}) Q^n_g, \label{mgre1-dis}\\
\f{C_v}{\mathcal{C}\mathcal{P}_0}\f{T^{n+1}-T^{n}}{\Delta t}   =  \sum\limits_{g'=1}^{G } \left[4\pi\mathscr{L}_a\left(\sigma^1_{a,g'}  \rho^{n+1}_{g'}- \sigma^2_{a,g'}  B^{n+1}_{g'}\right)\right],\label{mgre2-dis}
	\end{numcases}
\end{subequations}
Instead of solving  system \eqref{eqn:003re-dis} directly, we take the zeroth and first moments of \eqref{mgre1-dis}, and then couple them together with \eqref{mgre2-dis}. More precisely, 
\begin{subequations}\label{zero-one}  
	\begin{numcases}{}
	\f{1}{\mathcal{C}}\f{\rho_g^{n+1}-\rho_g^{n}}{\Delta t}+\average{\Omega \cdot \nabla ( \Omega\cdot R^{n+1}_g)}= \mathscr{L}_a\left(\sigma^2_{a,g}  B^{n+1}_g-\sigma^1_{a,g} \rho^{n+1}_g\right),  \label{zero-one--1}\\
	\f{1}{3\mathcal{C}}\f{R_g^{n+1}-R_g^{n}}{\Delta t}+\average{\Omega\Omega \cdot \nabla (\rho^{n+1}_g+Q^n_g)}= -\f{(\mathscr{L}_a+\mathscr{L}_s) \sigma_{t,g}R^{n+1}_g}{3
	},\label{zero-one--2}\\
\f{C_v}{\mathcal{C}\mathcal{P}_0}\f{T^{n+1}-T^{n}}{\Delta t}   =  \sum\limits_{g'=1}^{G } \left[4\pi\mathscr{L}_a\left(\sigma^1_{a,g'}  \rho^{n+1}_{g'}- \sigma^2_{a,g'}  B^{n+1}_{g'}\right)\right].\label{zero-one--3}
	\end{numcases}
\end{subequations}
The macroscopic quantities $\rho_g,\ R_g$ and $T$ can be updated by solving equation  \eqref{zero-one}, then $Q_g$ can be updated implicitly by equation \eqref{mgre1-dis}, i.e.,
\begin{equation}\label{eqnqg}
	\begin{aligned}
&	\f{1}{\mathcal{C}}\f{\rho_g^{n+1}-\rho_g^{n}}{\Delta t}+	\f{\Omega}{\mathcal{C}}\cdot\f{R_g^{n+1}-R_g^{n}}{\Delta t}+\f{1}{\mathcal{C}}\f{Q_g^{n+1}-Q_g^{n}}{\Delta t}+\Omega \cdot \nabla (\rho^{n+1}_g+\Omega\cdot R^{n+1}_g+Q^{n+1}_g )  =\\ &\hspace{1.0cm}\mathscr{L}_a\left(  \sigma^2_{a,g}B^{n+1}_g- \sigma^1_{a,g} \rho^{n+1}_g\right)-(\mathscr{L}_a+\mathscr{L}_s)\sigma_{t,g}\Omega\cdot R^{n+1}_g-(\mathscr{L}_a\sigma_{a,g}+\mathscr{L}_s\sigma_{s,g}) Q^{n+1}_g,
	\end{aligned}
\end{equation}
in which only a linear transport equation is solved with all direction and energy group being decoupled. Next, we will discuss the  diffusion limit of system \eqref{zero-one}.
\subsection{The gray radiation diffusion regime of  \eqref{zero-one}}
The asymptotic analysis of the time discretization is similar to the analysis in Section 3.3.  In the gray radiation diffusion regime, $\mathscr{L}_a=1/\eps$  and $\mathscr{L}_s=1/\eps$.
\textcolor{black}{
At the leading order of \eqref{zero-one--1}, one gets
\begin{equation}\label{rhogBgn--zj}
 \sigma^1_{a,g}\rho^{n+1}_g =  \sigma^2_{a,g} B^{n+1}_g + \order(\eps)\,,
\end{equation}
Summing the above equation over the group index $g$, and by the definitions \eqref{eqapp}, one can get
\begin{equation}
\int\limits_{0}^\infty\sigma_a\rho^{n+1}\text{d}\nu=\int\limits_{0}^\infty\sigma_aB^{n+1}\text{d}\nu + \order(\eps),
\end{equation}
which implies that
\begin{equation}\label{rhogBgn}
\rho^{n+1}_g = B^{n+1}_g + \order(\eps)\,.
\end{equation}}
In \eqref{zero-one--1}, adding up all the equations in each group and \eqref{zero-one--3}, we have
 \begin{equation} \label{diff-mg1re-semi-mi}
	4\pi 	\f{\sum_{g=1}^G \rho^{n+1}_g-\sum_{g=1}^G \rho^{n}_g}{\Delta t}+ C_v\f{T^{n+1}-T^{n}}{\Delta t}=  \f{4\pi}{3\eps}\sum_{g=1}^G\left(\nabla\cdot R_g^{n+1}\right),
\end{equation}
From \eqref{rhogBgn} and the conditions that $B_g$ satisfies in \eqref{Bgconditions}, $4\pi\sum_{g=1}^G \rho^{n+1}_g= (T^{n+1})^4+\order(\eps)$. To get the diffusion limit of equation \eqref{diff-mg1re-semi-mi}, the relation between $\rho_g^{n+1}$ and $R^{n+1}_g$ has to be analyzed.
By using Champan-Enskog expansion in equation \eqref{eqnqg},  one can get
\begin{align}
	- \f{\Omega\cdot R^{n+1}_g}{\int_{\nu_{g-\frac{1}{2}}}^{\nu_{g+\frac{1}{2}}}\f{1}{\sigma_{a}+ \sigma_{s}} \f{\partial B}{\partial T}(T^{n+1}) \text{d} \nu\big/\int_{\nu_{g-\frac{1}{2}}}^{\nu_{g+\frac{1}{2}}}\f{\partial B}{\partial T}(T^{n+1}) \text{d} \nu
	}-(\sigma_{a,g}+\sigma_{s,g}) Q^{n+1}_g=\order(\eps).\label{anatz1re-semi}
\end{align}
Taking the second moment in $\Omega$ of equation \eqref{anatz1re-semi} yields
\begin{equation}\label{kqg-semi}
	\average{\Omega\Omega Q^{n+1}_g}:=K^{n+1}_{Q_g}=\order(\eps).
\end{equation}
 Taking the first moment of \eqref{mgre1-dis}, and using the conditions $\rho^{n+1}_g\approx B^{n+1}_g$, $K^{n}_{Q_g}=\order(\eps)$ and $\average{\Omega Q^n_g}$=0, one has
\begin{equation}
	-\f{1}{\eps}R^{n+1}_g=\frac{1}{\int_{\nu_{g-\frac{1}{2}}}^{\nu_{g+\frac{1}{2}}}\f{\partial B}{\partial T} (T^{n+1})\text{d} \nu}\int_{\nu_{g-\frac{1}{2}}}^{\nu_{g+\frac{1}{2}}}\f{\nabla\rho^{n+1}_g}{\sigma_{a}+\sigma_{s} } \f{\partial B}{\partial T}(T^{n+1})\text{d} \nu+\order(\eps),
\end{equation}
which gives the relation between $R_g$ and $\rho_g$, then  putting it into \eqref{diff-mg1re-semi-mi},  and  using \eqref{rhogBgn}, when $\eps\to 0$, \eqref{diff-mg1re-semi-mi} reduces to
 \begin{equation} \label{diff-mg1re-semi}
   4\pi 	\f{\sum\limits_{g=1}^G B^{n+1}_g-\sum\limits_{g=1}^G B^{n}_g}{\Delta t}+ C_v\f{T^{n+1}-T^{n}}{\Delta t}=  \sum_{g=1}^G\nabla \cdot\left[ \left(\int\limits_{\nu_{g-1/2}}^{\nu_{g+1/2}}  \frac{4\pi}{3(\sigma_a+\sigma_s)} \f{\partial B}{\partial T}(T^{n+1})\text{d} \nu\right)\nabla T^{n+1}\right]  \,.
   \end{equation}
Simultaneously, we can obtain 
 \begin{equation} \label{diff-mg1re-semi01}
		\f{(T^{n+1})^4-(T^{n})^4}{\Delta t}+ C_v\f{T^{n+1}-T^{n}}{\Delta t}=  \nabla \cdot \left(  \frac{1}{3\widehat\sigma_R} \nabla (T^{n+1})^4\right) \,,
	\end{equation}
with
$$
\f{1}{\widehat\sigma_R}=\f{1}{4(T^{n+1})^3} \left(\sum\limits_{g=1}^{G }\int\limits_{\nu_{g-1/2}}^{\nu_{g+1/2}}  \frac{4\pi}{\sigma_a+\sigma_s} \f{\partial B}{\partial T}(T^{n+1})\text{d} \nu\right),
$$
which is a semi-discretization  for equation \eqref{diff-mg2}.
\subsection{The frequency dependent diffusion regime of  \eqref{zero-one}}
In the frequency dependent diffusion regime, $\mathscr{L}_a=\eps$  and $\mathscr{L}_s=1/\eps$. By using Champan-Enskog expansion in equation \eqref{eqnqg}, one can get
\begin{align}
	- \f{\Omega\cdot R^{n+1}_g}{\int_{\nu_{g-\frac{1}{2}}}^{\nu_{g+\frac{1}{2}}}\f{1}{\sigma_{s} }\f{\partial B}{\partial T}(T^{n+1}) \text{d} \nu\big/\int_{\nu_{g-\frac{1}{2}}}^{\nu_{g+\frac{1}{2}}}\f{\partial B}{\partial T}(T^{n+1}) \text{d} \nu
	}-\sigma_{s,g} Q^{n+1}_g=\order(\eps).\label{anatz1re-semin}
\end{align}
Taking the second moment for $\Omega$ in equation \eqref{anatz1re-semin}  yields
\begin{equation}\label{kqg-semin}
	\average{\Omega\Omega Q^{n+1}_g}=\order(\eps).
\end{equation}
From the second equation of system \eqref{zero-one}, the leading order is
\begin{equation}
	-\f{1}{\eps}R^{n+1}_g=\frac{1}{\int_{\nu_{g-\frac{1}{2}}}^{\nu_{g+\frac{1}{2}}}\f{\partial B}{\partial T}(T^{n+1}) \text{d} \nu}\int_{\nu_{g-\frac{1}{2}}}^{\nu_{g+\frac{1}{2}}}\f{\nabla\rho^{n+1}_g}{\sigma_{s} }\f{\partial B}{\partial T}(T^{n+1}) \text{d} \nu+\order(\eps),
\end{equation}
then putting it into the zeroth moment of \eqref{mgre1-dis}, and sending $\eps\to 0$, one gets
\begin{equation}
	\f{\rho_g^{n+1}-\rho_g^n}{\Delta t}-\nabla \cdot \left(\f{1}{3\sigma'_{s,g} }\nabla   \rho^{n+1}_g \right) =\sigma^2_{a,g}B_g^{n+1}-\sigma^1_{a,g}\rho_g^{n+1},
\end{equation}
where $\f{1}{\sigma'_{s,g}}=\int_{\nu_{g-\frac{1}{2}}}^{\nu_{g+\frac{1}{2}}}\f{1}{\sigma_{s} }\f{\partial B}{\partial T}(T_r^{n+1}) \text{d} \nu\big/\int_{\nu_{g-\frac{1}{2}}}^{\nu_{g+\frac{1}{2}}}\f{\partial B}{\partial T}(T_r^{n+1}) \text{d} \nu.$ Moreover, sending $\eps\to 0$ in equation \eqref{mgre2-dis} gives
\begin{equation}
	C_v\f{T^{n+1}-T^{n}}{\Delta t}   =  \sum\limits_{g'=1}^{G }\left[4\pi\left(  \sigma^1_{a,g'} \rho^{n+1}_{g'}- \sigma^2_{a,g'}  B^{n+1}_{g'}\right)\right].
\end{equation}
In summary, the diffusion limit of system \eqref{zero-one} in the frequency dependent diffusion regime is 
\begin{equation*}
	\left\{
	\begin{aligned}
		&		\f{\rho_g^{n+1}-\rho_g^n}{\Delta t}-\nabla \cdot \left(\f{1}{3\sigma'_{s,g} }\nabla   \rho^{n+1}_g \right) =\sigma^2_{a,g}B_g^{n+1}-\sigma^1_{a,g}\rho_g^{n+1}\,,\\
		& C_v\f{T^{n+1}-T^{n}}{\Delta t}   =  \sum\limits_{g'=1}^{G }\left[4\pi \left(  \sigma^1_{a,g'}\rho^{n+1}_{g'}-  \sigma^2_{a,g'}B^{n+1}_{g'}\right)\right],
	\end{aligned}
	\right.
\end{equation*}
which is a semi-discretization  for system \eqref{eqn:003re-nonequi-diff}.
\section{AP full discretization for FRTE}
For the ease of exposition, we will explain our spatial discretion in 1D. That is, $x \in [0,L]$, $\Omega \in [-1,1]$, and $\average{f(\Omega)} = \frac{1}{2} \int_{-1}^1 f(\Omega) \rd \Omega$, and the
 boundary conditions for the FRTE are
\begin{equation} \label{BC1d}
\begin{aligned}
I_g(t,0,\Omega) = b_{g,\text{L}}(t,\Omega)=\int\limits_{\nu_{g-1/2}}^{\nu_{g+1/2}}b_{\text{L}}(t,\Omega,\nu)\text{d}\nu, ~ \text{ for } ~\Omega>0; \\
I_g(t,L,\Omega) = b_{g,\text{R}}(t,\Omega)=\int\limits_{\nu_{g-1/2}}^{\nu_{g+1/2}}b_{\text{R}}(t,\Omega,\nu)\text{d}\nu, ~ \text{ for } ~ \Omega<0\,.
\end{aligned}
\end{equation}
Higher dimensions can be treated in the dimension by dimension manner. 
Let $\Delta x=L/N_x$, and we consider the uniform mesh as follows
 \[
 x_{i}=(i-1)\Delta x, \quad i = 1, 2 \cdots, N_{x}+1,
 \]
 and let
 $$x_{i+\f{1}{2}}=\left(x_{i}+x_{i+1}\right)/2,\qquad \mbox{for $i=1,\cdots, N_x.$}$$
\subsection{Even-odd space discretization}To get a consistent stencil in spatial discretization, we use  the even-odd parity method \cite{k2016asymptotic,tang2021accurate,jin2000uniformly}. Let
\begin{equation} \label{EO}
	E_g(\Omega)=\f{1}{2}\big(I_g(\Omega)+I_g(-\Omega)\big),  \qquad O_g(\Omega)=\frac{1}{2} (I_g(\Omega)-I_g(-\Omega)),\qquad  \Omega>0
\end{equation}
be respectively the even and odd part of $I_g$, and we consider the even part on half spatial grid and odd part on regular grid, i.e., 
	\begin{align}
		& E_{g,i+1/2} (\Omega)\approx E_g(x_{i+1/2}, \Omega), \quad i = 1, \cdots, N_x\,; \label{E000}
		\\ & O_{g,i} (\Omega)\approx O_g(x_i,\Omega), \quad i = 1, \cdots, N_x + 1\,. \label{O000}
	\end{align}
Moreover, plugging the decomposition \eqref{eqn:011} into \eqref{EO}, we have 
\begin{equation}\label{eqn::eqoq}
E_g = \rho_g + \big(Q_g(\Omega)+Q_g(-\Omega)\big) := \rho_g+  E_{Q_g}, \quad O_g =\Omega\cdot R_g+ \f{1}{2}\big( Q_g(\Omega)-Q_g(-\Omega)\big) :=\Omega\cdot R_g+O_{Q_g}\,.
\end{equation}
Then the full discretization for the macroscopic equation \eqref{zero-one}   and the microscopic equation \eqref{eqnqg} are described below.

\vspace{0.5cm}

\textbf{Space discretization for the first term \eqref{zero-one}.} System \eqref{eqn:003re-dis} can be written into a system for $\rho_g,\ R_g$ and $E_{Q_g},\ O_{Q_g}$:
\begin{equation} \label{9163}
	\left\{
	\begin{aligned}
		& \frac{1}{\mathcal{C}} \partial_t (\Omega\cdot R_g+O_{Q_g})  + \Omega\cdot \nabla_x (\rho_g+  E_{Q_g}) = -(\mathscr{L}_a+\mathscr{L}_s)\sigma_{t,g}\Omega\cdot R_g-(\mathscr{L}_a\sigma_{a,g}+\mathscr{L}_s\sigma_{s,g}) O_{Q_g}\,;
		\\ & \frac{1}{\mathcal{C}} \partial_t (\rho_g+  E_{Q_g})+ \Omega \cdot \nabla_x (\Omega\cdot R_g+O_{Q_g})= \mathscr{L}_a \left(  \sigma^2_{a,g}B_g-  \sigma^1_{a,g}\rho_g\right)-(\mathscr{L}_a\sigma_{a,g}+\mathscr{L}_s\sigma_{s,g}) E_{Q_g}\,;
		\\ &  	\f{C_v}{\mathcal{C}\mathcal{P}_0}\partial_t T = \sum\limits_{g'=1}^{G } \left[4\pi\mathscr{L}_a \left( \sigma^1_{a,g'}\rho_{g'}- \sigma^2_{a,g'} B_{g'}\right)\right]\,.
	\end{aligned}
	\right.
\end{equation}
In terms of $\rho_g,\ R_g$ and $E_{Q_g},\ O_{Q_g}$, the system \eqref{zero-one}  becomes
\begin{equation} \label{9163-m}
	\left\{
\begin{aligned}
	&	\f{1}{3\mathcal{C}}\partial_t R_g+\average{\Omega\Omega \cdot \nabla (\rho_g+E_{Q_g})}= -\f{(\mathscr{L}_a+\mathscr{L}_s)\sigma_{t,g} R_g}{3
	},\\
&	\f{1}{\mathcal{C}}\partial_t \rho_g+\average{\Omega \cdot \nabla ( \Omega\cdot R_g)}= \mathscr{L}_a\left(\sigma^2_{a,g}  B_g- \sigma^1_{a,g}\rho_g\right), \\
	&\f{C_v}{\mathcal{C}\mathcal{P}_0}\partial_t  T =  \sum\limits_{g'=1}^{G } \left[4\pi\mathscr{L}_a \left( \sigma^1_{a,g'}\rho_{g'}- \sigma^2_{a,g'} B_{g'}\right)\right].
\end{aligned}
\right.
\end{equation}
The spatial discretization then takes the following form: 
\begin{subequations} \label{dis-mac}
	\begin{numcases}{}
			\f{1}{3\mathcal{C}} \frac{R_{g,i}^{n+1}-R_{g,i}^{n}}{\Delta t}+\frac{1}{3\Delta x}(\rho^{n+1}_{g,i+1/2}-\rho^{n+1}_{g,i-1/2})+\frac{1}{\Delta x}\left(\int_0^1\Omega\Omega E^n_{Q_g,i+1/2}\text{d}\Omega-\int_0^1\Omega\Omega E^n_{Q_g,i-1/2}\text{d}\Omega\right)\nonumber \\
				\hspace{7.5cm}= -\f{(\mathscr{L}_a+\mathscr{L}_s)\sigma_{t,g,i}R^{n+1}_{g,i}}{3},\quad 2 \leq i \leq N_x \,; \label{dis-cor1}\\
	\f{1}{\mathcal{C}} \frac{\rho_{g,i+1/2}^{n+1}-\rho_{g,i+1/2}^{n}}{\Delta t}+ \frac{R^{n+1}_{g,i+1}-R^{n+1}_{g,i}}{3\Delta x}  =\mathscr{L}_a\left(\sigma^2_{a,g,i+1/2}   B^{n+1}_{g,i+1/2}-\sigma^1_{a,g,i+1/2}  \rho^{n+1}_{g,i+1/2}\right),\nonumber \\
	\hspace{12cm} \quad 1 \leq i \leq N_x\,; \label{dis-cor2}\\
	\f{C_v}{\mathcal{C}\mathcal{P}_0}\frac{ T_{i+1/2}^{n+1}-T_{i+1/2}^n}{\Delta t}=\sum\limits_{g'=1}^{G } \left[4\pi\mathscr{L}_a \left(\sigma^1_{a,g',i+1/2} \rho^{n+1}_{g',i+1/2}-  \sigma^2_{a,g',i+1/2}B^{n+1}_{g',i+1/2}\right)\right],  1\leq i \leq N_x \,. \label{dis-cor3}
	\end{numcases}
\end{subequations}
\vspace{0.5cm}

\textbf{Space discretization for the second term \eqref{eqnqg}.}  The variables $\rho_g^{n+1}$, $R_g^{n+1}$ and $T^{n+1}$ will be updated by solving system \eqref{dis-mac}, and this is a nonlinear system that has to be solved by Newton iteration. Then the full discretization of \eqref{eqnqg} to update $O_{Q_g}$, $E_{Q_g}$ is:
\begin{subequations} \label{eqn:upq}
\begin{numcases}{}
	\frac{1}{\mathcal{C}} \f{O^{n+1}_{Q_g,i}- O^{n}_{Q_g,i}}{\Delta t} +\frac{\Omega}{\mathcal{C}}\f{R^{n+1}_{g,i}- R^{n}_{g,i}}{\Delta t}+ \Omega\f{\rho^{n+1}_{g,i+1/2}-\rho^{n+1}_{g,i-1/2}}{\Delta x} + \Omega\f{E^{n+1}_{Q_g,i+1/2}-E^{n+1}_{Q_g,i-1/2}}{\Delta x}  \nonumber \\
		\hspace{1.5cm}= -(\mathscr{L}_a+\mathscr{L}_s)\sigma_{t,g,i}\Omega\cdot R^{n+1}_{g,i}-(\mathscr{L}_a\sigma_{a,g,i}+\mathscr{L}_s\sigma_{s,g,i}) O^{n+1}_{Q_g,i},\quad 2 \leq i \leq N_x \,; \label{dis-cor1-mic}\\
		\f{1}{\mathcal{C}} \frac{\rho_{g,i+1/2}^{n+1}-\rho_{g,i+1/2}^{n}}{\Delta t}+\frac{1}{\mathcal{C}} \f{E^{n+1}_{Q_g,i+1/2}- E^{n}_{Q_g,i+1/2}}{\Delta t}+ \Omega\Omega\f{R^{n+1}_{g,i+1}-R^{n+1}_{g,i}}{\Delta x} + \Omega\f{O^{n+1}_{Q_g,i+1}-O^{n+1}_{Q_g,i}}{\Delta x}  \nonumber \\ 	\hspace{1cm} =\mathscr{L}_a \left( \sigma^2_{a,g,i+1/2} B^{n+1}_{g,i+1/2}- \sigma^1_{a,g,i+1/2}\rho^{n+1}_{g,i+1/2}\right)-(\mathscr{L}_a\sigma_{a,g,i+1/2}+\mathscr{L}_s\sigma_{s,g,i+1/2}) E^{n+1}_{Q_g,i+1/2},\nonumber \\ 	\hspace{11cm}
	\quad 1 \leq i \leq N_x\,. \label{dis-cor2-mic}
	\end{numcases}
\end{subequations}
Here only  a linear transport equation is solved with all direction and energy group being decoupled.

\textbf{Boundary conditions}  To cope with \eqref{BC1d}, using the relation \eqref{EO}, we have 
\begin{equation} 
E_{g,3/2}^{n+1} + \frac{1}{2} \left( O_{g,2}^{n+1}  + O_{g,1}^{n+1} \right) = b_\text{g,L}(\Omega), ~ \Omega>0; \quad
E_{g,N_x+1/2}^{n+1} - \frac{1}{2} \left( O_{g,N_x}^{n+1}  + O_{g,N_x+1}^{n+1} \right) = b_\text{g,R}(-\Omega), ~ \Omega>0\,.
\end{equation}
Hence,
\begin{equation} \label{bc-pre}
O_{g,1}^{n+1}(\Omega) = 2 \left( b_\text{g,L} (\Omega) - E_{g,3/2} ^{n+1}\right) - O_{g,2}^{n+1}, \quad 
O_{g,N_x+1}^{n+1}(\Omega) = 2 \left( E_{g,N_x + \half} ^{n+1}-b_\text{g,R} (-\Omega) \right) - O_{g,N_x}^{n+1}.
\end{equation}
To summarize, we have the following one time step update of the fully discrete version of FRTE.
\begin{algorithm}[!h]
\caption{one step of fully discrete update for FRTE}
\SetAlgoLined
\KwIn{$T_{i+1/2}^n$, $E_{g,i+1/2}^s$,$O_{g,i}^n$ }
\KwOut{$T_{i+1/2}^{n+1}$,$E_{g,i+1/2}^{n+1}$,$O_{g,i}^{n+1}$}
\BlankLine
      get $R_{g,i}^{n+1}$, $\rho_{g,i+1/2}^{n+1}$, $T_{ i+1/2}^{n+1}$ from the system \eqref{dis-mac};
    \\  get $O_{Q_g,i}^{n+1}$ and $E_{Q_g,i+1/2}^{n+1}$ from the system   \eqref{eqn:upq};
    \\  obtain $E_{g,i+1/2}^{n+1}$,$O_{g,i}^{n+1}$ from \eqref{eqn::eqoq}.
\end{algorithm}
\subsection{The gray radiation diffusion regime of \eqref{dis-mac}} In this subsection, the gray radiation diffusion regime is considered, i.e., $\mathscr{L}_a=1/\eps$  and $\mathscr{L}_s=1/\eps$. 
From \eqref{dis-cor2}, one can get
\begin{equation}\label{r-d-dis--zj}
 \sigma^1_{a,g,i+1/2}\rho^{n+1}_{g,i+1/2}= \sigma^2_{a,g,i+1/2}B^{n+1}_{g,i+1/2}+\order(\eps),
\end{equation}
Summing the above equation over the group index $g$, and by the definitions \eqref{eqapp}, one can get
\begin{equation}
\int\limits_{0}^\infty\sigma_a\rho_{g,i+1/2}^{n+1}\text{d}\nu=\int\limits_{0}^\infty\sigma_aB_{g,i+1/2}^{n+1}\text{d}\nu + \order(\eps),
\end{equation}
which implies that
\begin{equation}\label{r-d-dis}
\rho^{n+1}_{g,i+1/2}= B^{n+1}_{g,i+1/2}+\order(\eps).
\end{equation}
 Next, we will find the relation between $\rho^{n+1}_{g,i+1/2}$ and $R^{n+1}_{g,i}$.  By using Champan-Enskog expansion in equation \eqref{dis-cor2-mic}, and using \eqref{r-d-dis}, one can get
\begin{align}
	-\sigma_{a,g,i+1/2} E_{Q_g,i+1/2}^{n+1}=\order(\eps).\label{anatz1re-full}
\end{align}
Taking the second moment for $\Omega$ in equation \eqref{anatz1re-full}  yields
\begin{equation}\label{kqg-full}
	\average{\Omega\Omega E_{Q_g,i+1/2}^{n+1}}=\order(\eps).
\end{equation}
By using \eqref{kqg-full}, and letting $\eps\to0$, equation \eqref{dis-cor1} reduces to
\begin{equation}\label{r-d-dis1}
-\f{1}{\eps}R^{n+1}_{g,i}=\frac{1}{\int_{\nu_{g-\frac{1}{2}}}^{\nu_{g+\frac{1}{2}}} \f{\partial B}{\partial T}(T_i^{n+1})\text{d} \nu}\int_{\nu_{g-\frac{1}{2}}}^{\nu_{g+\frac{1}{2}}}\f{1}{\sigma_{a,i}+\sigma_{s,i} }\f{\partial B}{\partial T}(T_i^{n+1}) \text{d} \nu\f{\rho^{n+1}_{g,i+1/2}-\rho^{n+1}_{g,i-1/2}}{\Dx}.
\end{equation}
Substituting \eqref{r-d-dis1} into equation \eqref{dis-cor2}, adding up all the equations in each group, then adding it to  equation  \eqref{dis-cor3}, using the equation \eqref{r-d-dis} and letting $\eps\to 0$, one can get
\begin{equation}
	\begin{aligned}
     &4\pi 	\f{\sum_{g=1}^G B^{n+1}_{g,i+1/2}-\sum_{g=1}^G B^{n}_{g,i+1/2}}{\Delta t}+ C_v\f{T_{i+1/2}^{n+1}-T_{i+1/2}^{n}}{\Delta t} =  \sum_{g=1}^G\frac{1}{3\Delta x}\\
    &\left[\f{\int_{\nu_{g-\frac{1}{2}}}^{\nu_{g+\frac{1}{2}}}\f{4\pi}{\sigma_{a,i+1}+\sigma_{s,i+1} }\f{\partial B}{\partial T}(T_{i+1}^{n+1})\text{d} \nu}{\int_{\nu_{g-\frac{1}{2}}}^{\nu_{g+\frac{1}{2}}} \f{\partial B}{\partial T}(T_{i+1}^{n+1})\text{d} \nu} \f{B^{n+1}_{g,i+3/2}-B^{n+1}_{g,i+1/2}}{\Dx}-\f{\int_{\nu_{g-\frac{1}{2}}}^{\nu_{g+\frac{1}{2}}}\f{4\pi}{\sigma_{a,i}+\sigma_{s,i} }\f{\partial B}{\partial T}(T_{i}^{n+1})\text{d} \nu}{\int_{\nu_{g-\frac{1}{2}}}^{\nu_{g+\frac{1}{2}}} \f{\partial B}{\partial T}(T_{i}^{n+1})\text{d} \nu} \f{B^{n+1}_{g,i+1/2}-B^{n+1}_{g,i-1/2}}{\Dx}\right] \,.
 \end{aligned}
\end{equation}
Moreover, we have
$$
4\pi \sum_{g=1}^G B^{n+1}_{g,i+1/2}=\int_{0}^{\infty}B(\nu,T^{n+1}_{i+1/2})\text{d}\nu=(T^{n+1}_{i+1/2})^4,
$$
which is is consistent with equation \eqref{diff-mg2}, and the calculations of the other two gray radiation diffusion states are similar.
\subsection{The frequency dependent diffusion regime of \eqref{dis-mac}} In the frequency dependent diffusion regime, $\mathscr{L}_a=\eps$  and $\mathscr{L}_s=1/\eps$. 
By using Champan-Enskog expansion in equation \eqref{dis-cor2-mic}, one can get
\begin{align}
	-\sigma_{s,g,i+1/2} E_{Q_g,i+1/2}^{n+1}=\order(\eps).\label{anatz1re-full-n}
\end{align}
Taking the second moment for $\Omega$ in equation \eqref{anatz1re-full-n} yields
\begin{equation}\label{kqg-full-n}
	\average{\Omega\Omega E_{Q_g,i+1/2}^{n+1}}=\order(\eps).
\end{equation}
By using \eqref{kqg-full-n}, and letting $\eps\to0$, equation \eqref{dis-cor1} reduces to
\begin{equation}\label{r-d-dis1-n}
	-\f{1}{\eps}R^{n+1}_{g,i}=\frac{1}{\int_{\nu_{g-\frac{1}{2}}}^{\nu_{g+\frac{1}{2}}}\f{\partial B}{\partial T}(T_{r,i}^{n+1}) \text{d} \nu}\int_{\nu_{g-\frac{1}{2}}}^{\nu_{g+\frac{1}{2}}}\f{1}{\sigma_{s,i} }\f{\partial B}{\partial T}(T_{r,i}^{n+1}) \text{d} \nu\f{\rho^{n+1}_{g,i+1/2}-\rho^{n+1}_{g,i-1/2}}{\Dx}.
\end{equation}
Substituting \eqref{r-d-dis1-n} into equation \eqref{dis-cor2}, and letting $\eps\to 0$, one can get
\begin{equation}
	\begin{aligned}	
	\f{\rho_{g,i+1/2}^{n+1}-\rho_{g,i+1/2}^{n}}{\Delta t}&-\f{1}{3\Dx}\left[\f{\int_{\nu_{g-\frac{1}{2}}}^{\nu_{g+\frac{1}{2}}}\f{1}{\sigma_{s,i+1}}\f{\partial B}{\partial T}(T_{r,i+1}^{n+1})\text{d} \nu}{\int_{\nu_{g-\frac{1}{2}}}^{\nu_{g+\frac{1}{2}}}\f{\partial B}{\partial T}(T_{r,i+1}^{n+1}) \text{d} \nu}\f{\rho^{n+1}_{g,i+3/2}-\rho^{n+1}_{g,i+1/2}}{\Dx}\right.\\
	&\hspace{0cm}\left.-\f{\int_{\nu_{g-\frac{1}{2}}}^{\nu_{g+\frac{1}{2}}}\f{1}{\sigma_{s,i}}\f{\partial B}{\partial T}(T_{r,i}^{n+1})\text{d} \nu}{\int_{\nu_{g-\frac{1}{2}}}^{\nu_{g+\frac{1}{2}}}\f{\partial B}{\partial T}(T_{r,i}^{n+1}) \text{d} \nu}\f{\rho^{n+1}_{g,i+1/2}-\rho^{n+1}_{g,i-1/2}}{\Dx}\right]= \sigma^2_{a,g,i+\f{1}{2}} B^{n+1}_{g,i+\f{1}{2}}- \sigma^1_{a,g,i+\f{1}{2}} \rho^{n+1}_{g,i+\f{1}{2}}.
		\end{aligned}
\end{equation}
Moreover, letting $\eps\to0$, equation \eqref{dis-cor3} reduces to
\begin{equation}
	C_v\f{T_{i+1/2}^{n+1}-T_{i+1/2}^{n}}{\Delta t}   =  \sum\limits_{g'=1}^{G }4\pi \left(\sigma^1_{a,g',i+\f{1}{2}} \rho^{n+1}_{g',i+\f{1}{2}}- \sigma^2_{a,g',i+\f{1}{2}} B^{n+1}_{g',i+\f{1}{2}}\right).
\end{equation}
which is consistent with equation \eqref{eqn:003re-nonequi-diff}.
\section{Numerical Examples}
In this section, we conduct several numerical examples to test the performance of the AP scheme for the FRTE. The units of the length, time, temperature and energy are respectively $cm,\ ns, \ keV$ and $Jk$ ($1\ Jk=10^9\ J$). With the above units, let the speed of light $c$ be 29.98\ $cm\ ns^{-1}$, the specific heat capacity $C_v$ be 0.1\ $Jk\ keV^{-1}\ cm^{-3}$ and the radiation constant $a_r$ be 0.01372\ $Jk\ keV^{-4}\ cm^{-3}$. As in \cite{densmore2012hybrid,sunj}, in all of these examples, the opacities take the following forms :
\begin{equation}\label{eqnsa}
	\sigma_{a}=\f{\sigma_{a0}(x)}{(h\nu)^3\sqrt{kT}},\qquad \sigma_{s}=\f{\sigma_{s0}(x)}{(h\nu)^3\sqrt{kT}}.
\end{equation}
In all numerical examples in this section, similar as in \cite{densmore2012hybrid,sunj}, we divide the frequency domain from 0.1 $eV$ to 100 $keV$ into 30 frequency groups, and the groups are logarithmically spaced. In each group,  we determine the group absorption, scattering and total cross section by using \eqref{eq:groupcoefficient} and \eqref{eqnsa}. Moreover, from Figure \ref{figequ-2sa}, the piece-wise constant approximations $\sigma_{a,g}^c$ are closer to the reference solutions, thus we use $\omega_1=\omega_2=1$ in all examples.
\subsection{Example 1}
In this example, we will show the convergence order and stability of the proposed scheme.
The thickness of the computational domain is $2\ cm$ and
the initial temperature is
\begin{equation} \label{IC000}
	T(0,x)=\max\big\{1-20(x-1)^2,10^{-3}\big\}\ keV.
\end{equation}
The initial radiation intensity is set to be a Planck distribution evaluated at the  temperature \eqref{IC000}. Zero Dirichlet boundary condition is used in this example. Four different sets of $\sigma_{a0}$, $\sigma_{s0}$ are tested: 1) $\sigma_{a0}=1\ keV^{7/2}\ cm^{-1}$, $\sigma_{s0}=1000\ keV^{7/2}\ cm^{-1}$; 2) $\sigma_{a0}=1000\ keV^{7/2}\ cm^{-1}$, $\sigma_{s0}=1\ keV^{7/2}\ cm^{-1}$; 3) $\sigma_{a0}=1000\ keV^{7/2}\ cm^{-1}$, $\sigma_{s0}=1000\ keV^{7/2}\ cm^{-1}$;
4) $\sigma_{a0}=1\ keV^{7/2}\ cm^{-1}$, $\sigma_{s0}=1\ keV^{7/2}\ cm^{-1}$. In the four cases, the first case belongs to the frequency dependent diffusion regime, and the second case and the third case are the gray radiation diffusion regime, and $\sigma_a$ and $\sigma_s$ are both optical thin in the last case.

\begin{figure}[ht]
	\centering
	\includegraphics[width=0.7\textwidth]{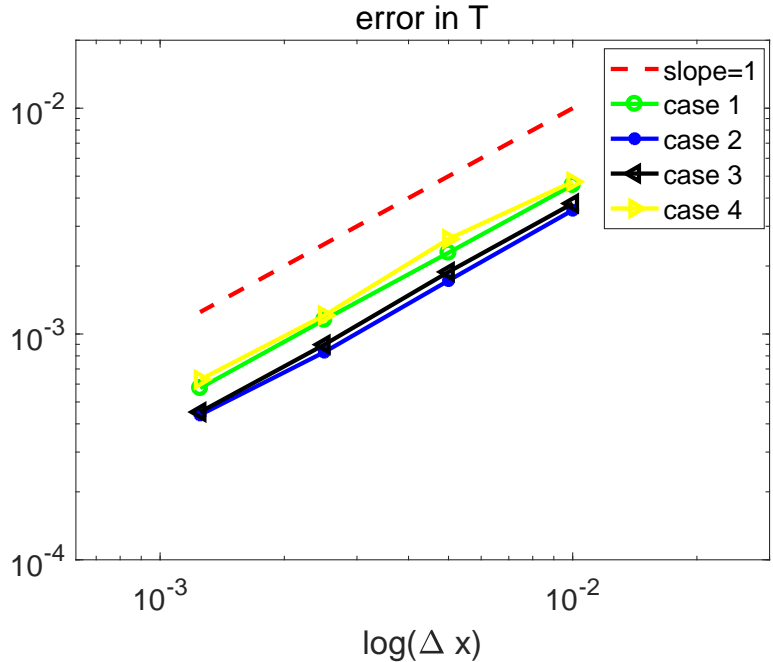}
		\caption{Example 1. Errors of different $\Delta x$, $\Delta t$ for the four cases. Here different $\Delta x=2*10^{-2}\ cm$, $10^{-2}\ cm$, $5*10^{-3}\ cm$, $2.5*10^{-3}\ cm$, $1.25*10^{-3}\ cm$  are tested and  $\Delta t$ are chosen to be $\Delta t=2*10^{-2}\ ns$, $10^{-2}\ ns$, $5*10^{-3}\ ns$, $2.5*10^{-3}\ ns$, $1.25*10^{-3}\ ns$ respectively. }
		\label{fig:er1-case}
	\end{figure}

\begin{figure}[!h]
	\centering
	\subfigure[$\sigma_{a0}=1\ keV^{7/2}\ cm^{-1},\ \sigma_{s0}=1000\ keV^{7/2}\ cm^{-1}$]{
		\includegraphics[width=3.2in]{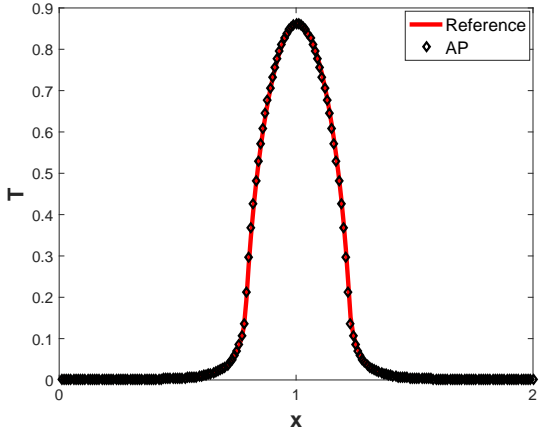}\label{fig.2_ab}
	}
	\subfigure[$\sigma_{a0}=1000\ keV^{7/2}\ cm^{-1},\ \sigma_{s0}=1\ keV^{7/2}\ cm^{-1}$]{
		\includegraphics[width=3.2in]{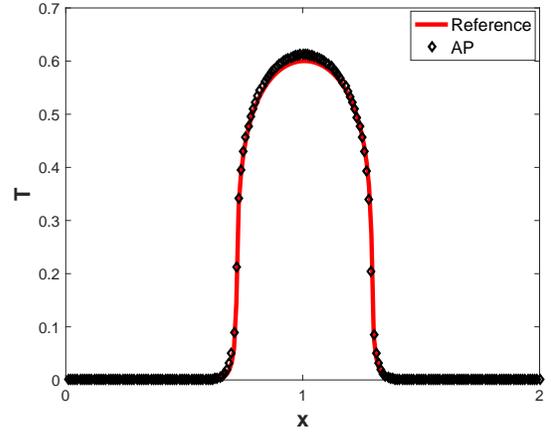}\label{fig.2_rt}
	}\\
	\subfigure[$\sigma_{a0}=1000\ keV^{7/2}\ cm^{-1}, \sigma_{s0}=1000\ keV^{7/2}\ cm^{-1}$]{
		\includegraphics[width=3.2in]{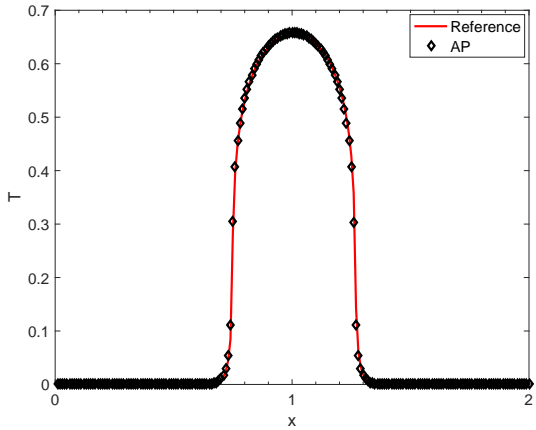}\label{fig.2_q1}
	}
	\subfigure[$\sigma_{a0}=1\ keV^{7/2}\ cm^{-1}, \sigma_{s0}=1\ keV^{7/2}\ cm^{-1}$]{
	\includegraphics[width=3.2in]{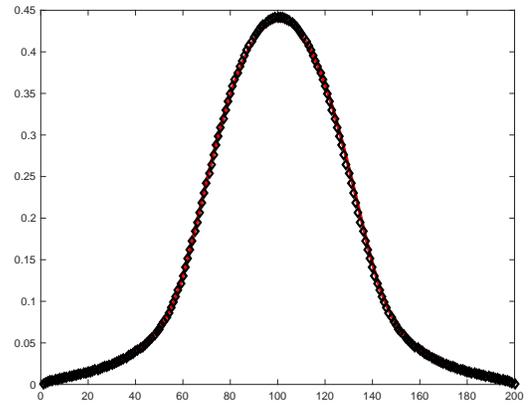}\label{fig.2_qq1}
}
	\caption{Example 1. Comparison of the material temperature using our AP scheme for equilibrium regime and the reference solution at $t=2\ ns$. }
		\label{figequ-2}
\end{figure}

In Fig.~\ref{fig:er1-case}, we plot the errors 
\begin{equation} \label{error}
	\text{error}_T = \| T_{\Delta x}(\cdot, t_\text{max}) - T_{\Delta x/2}(\cdot, t_\text{max})\|_{l_1}
\end{equation}
with $\Delta x=2*10^{-2}\ cm$, $10^{-2}\ cm$, $5*10^{-3}\ cm$, $2.5*10^{-3}\ cm$, $1.25*10^{-3}\ cm$ and $\Delta t=2*10^{-2}\ ns$, $10^{-2}\ ns$, $5*10^{-3}\ ns$, $2.5*10^{-3}\ ns$, $1.25*10^{-3}\ ns$, and $t_\text{max} = 1\ ns$. A uniform first order accuracy can be observed for all four different cases. Moreover, it is important to point out that the scheme is stable when the time step $\Delta t=8\Delta x$ for all the cases.
At last, we compare the results in the  gray radiation diffusion regime and frequency dependent diffusion regime to the reference results. 
For the four cases, we use a uniform spatial  mesh  and the cell size $\Delta x=0.01\ cm$, and the time step $\Delta t=0.04\ ns$. Moreover, the reference solution of all the cases are computed with $\Delta x=0.001\ cm$ and  $\Delta t=0.00001\ ns$. The numerical results at $t = 2\ ns$ are displayed in Fig.~\ref{figequ-2}.   
The material temperature from our AP method agrees well with the reference solution in all the four cases.

\subsection{Example 2}
\begin{figure}[!h]
	\centering
	\subfigure[$\sigma_{a0}=10\ keV^{7/2}\ cm^{-1}$]{
		\includegraphics[width=2.9in]{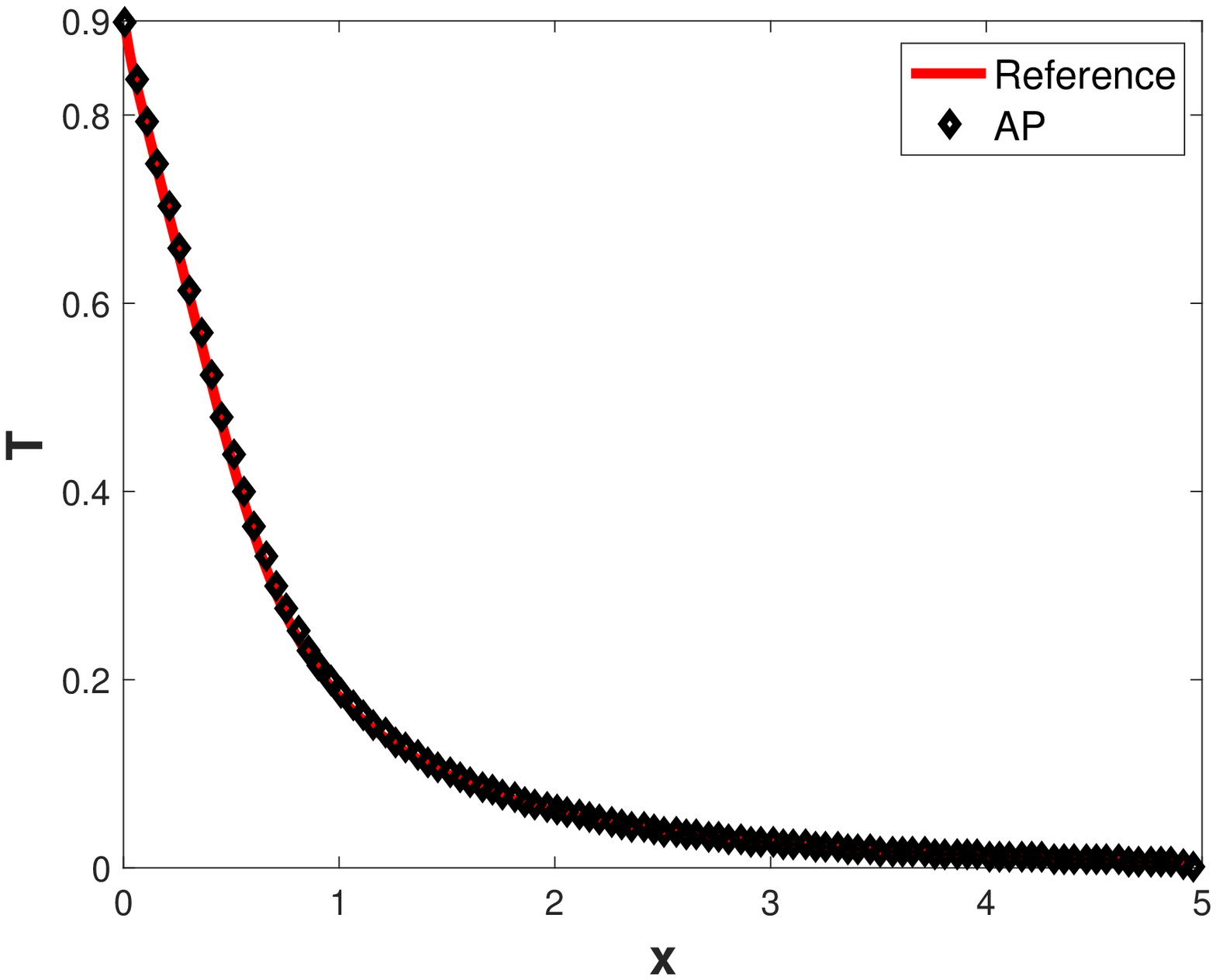}\label{fig.2_absa}
	}
	\subfigure[$\sigma_{a0}=100\ keV^{7/2}\ cm^{-1}$]{
		\includegraphics[width=2.9in]{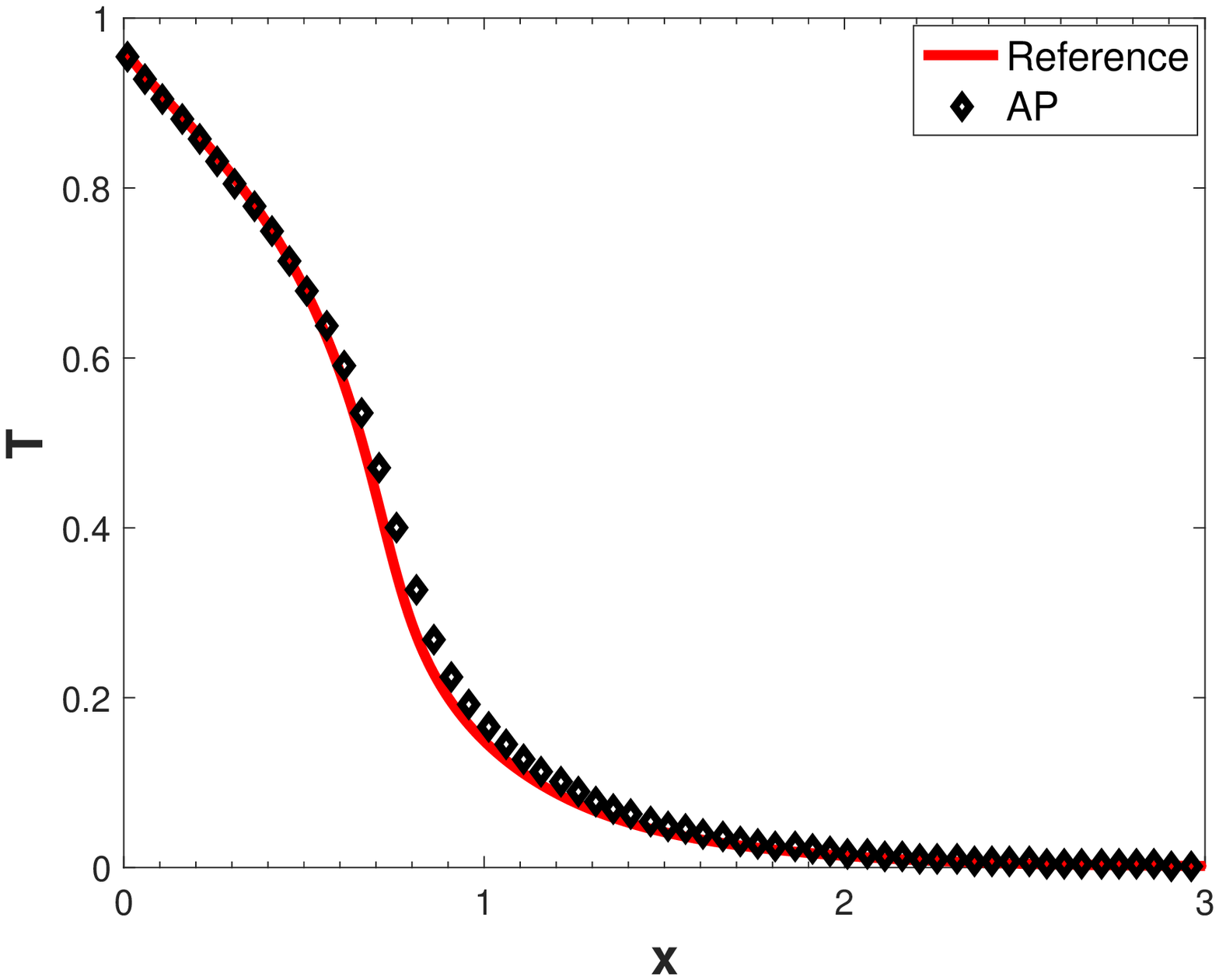}\label{fig.2_rtsa}
	}\\
	\subfigure[$\sigma_{a0}=1000\ keV^{7/2}\ cm^{-1}$]{
		\includegraphics[width=2.9in]{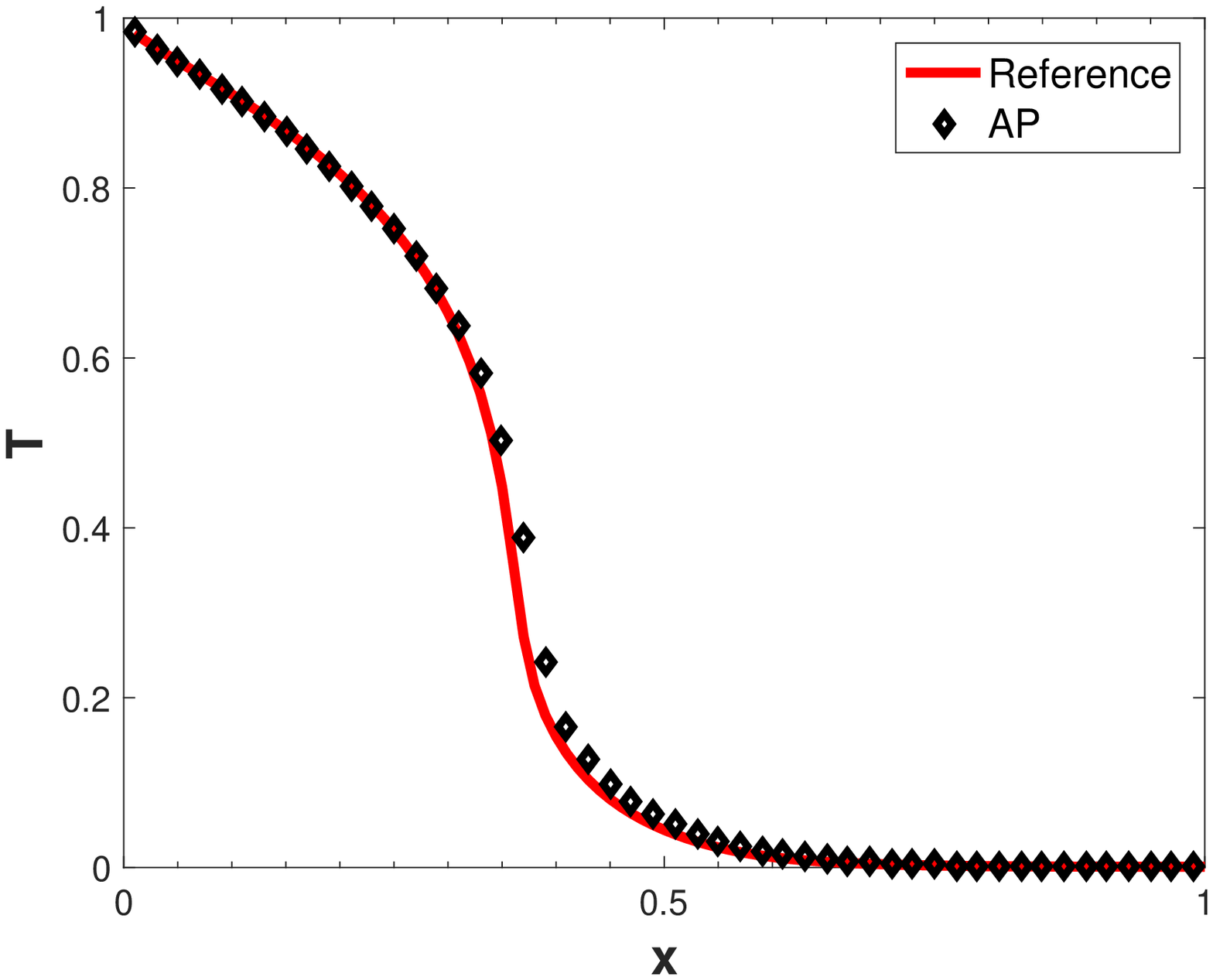}\label{fig.2_q1sa}
	}
	\caption{Example 2. Plot of three different cases at $t=1\ ns$. The red curves are reference solutions computed  in the finer mesh. Here $\Delta x=0.005\ cm$, in our AP scheme $\Delta t=0.005 \ ns$,  and in the reference solution, $\Delta t = 0.00001 \ ns$.}
	\label{fig:ho}
\end{figure}

We test an example in \cite{densmore2012hybrid,sunj,steinberg2022multi} where particles are injected from the left boundary. Three space homogeneous problems with $\sigma_{a0}=10\ keV^{7/2}\ cm^{-1},\ 100\ keV^{7/2}\ cm^{-1}$, $1000\ keV^{7/2}\ cm^{-1}$, and $\sigma_{s0}=0\ keV^{7/2}\ cm^{-1}$ are tested. The thickness of the computational domain is $5\ cm$,
the initial temperature is $10^{-3}\ keV$
and the initial radiation intensity is a Planck distribution evaluated at the initial temperature. On the left boundary, the incident radiation intensity  is given by a Planck distribution with the temperature being 1\ $keV$, and  we use a reflective boundary condition on the right boundary. For all three tests, we use a uniform spatial mesh with $\Delta x=0.005\ cm$ and $\Delta t=0.005 \ $ns to ensure the scheme  stable. The reference solution is calculated using a finer temporal mesh such that $\Delta t = 0.00001 \ $ns. Moreover,  in the example , the time size $\Delta t =0.005$ \ ns used in our AP method is much bigger than the time step $\Delta t =0.004/c$\ ns used in \cite{sunj}.
The numerical results at $t = 1\ ns$ are displayed in Fig.~\ref{fig:ho} (Note that we only show the part where the material temperature differs significantly from the initial value in the second and third figure).   
The good agreement between our solution and the reference solution indicates that our method works well in the optical thin and thick regime.
\subsection{Example 3}
\begin{figure}[ht]
	\centering
	\includegraphics[width=0.7\textwidth]{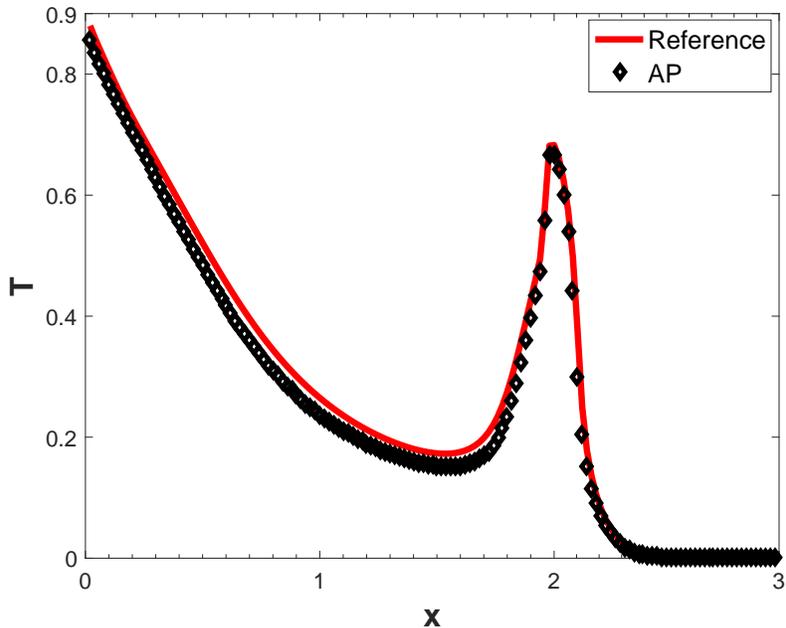}
	\caption{Example 3. Plot of multi-region problem at $t=1\ ns$. The red curves are reference solutions computed  in the finer mesh. Here $\Delta x=0.02\ cm$,  $\Delta t=0.02 \ ns$ in our AP scheme,  and in the reference solution, $\Delta x=0.005\ cm$ and   $\Delta t = 0.00001 \ ns$. }
	\label{fig:larsen}
\end{figure}

A 1D Marshak wave problem similar as in \cite{steinberg2022multi,densmore2012hybrid,sunj} is tested, where both optical thin and thick regions coexist. The opacities take  the forms of \eqref{eqnsa} with 
$$
\sigma_{a0}= \begin{cases}10 \quad\mathrm{keV}^{7/2}\ cm^{-1}, & 0<x<2.0 \ cm, \\ 1000\quad \mathrm{keV}^{7/2}\ cm^{-1}, & 2.0<x<3.0 ~cm,\end{cases}
$$
and $\sigma_{s0}=0$. Moreover, the width of the computational domain is $3$cm and we divide it uniformly into $150$ spatial cells. The time step is $\Delta t=0.02\ ns$.  The problem begins with flat material and radiation temperature profiles at 1 $eV$. A boundary source of 1 $keV$ is placed at the left boundary, while the right boundary is vacuum.  In this example, the reference solution is computed with $\Delta x=0.001\ cm$ and  $\Delta t=0.00001\ ns$. The numerical results at $t = 1\ ns$ are displayed in Fig.~\ref{fig:larsen}.  Moreover, the time size $\Delta t =0.02$ \ ns of our AP method in the example  is much bigger than the time step $\Delta t =0.016/c$\ ns used in \cite{sunj}. 
The material temperature from our AP method agrees well with the reference solution in the co-exit regime.

\section{Conclusion}
In this paper, we propose an AP decomposed multi-group scheme for the FRTE in the gray radiation diffusion regime and the frequency-dependent diffusion regime. The idea is to decompose the intensity into three parts: the zeroth moment of the intensity, the first moment of the intensity, and the residual term. The main contribution is the new decomposed multi-group energy discretization, which allows tuning of the approximations of the coefficients in the radiation diffusion regime. 

In the full AP discretization, the two macroscopic moments and the material temperature are updated first by treating them implicitly with the explicit residual term. After the zeroth and first-order moments are obtained, we update the residual term by solving implicitly a linear transport equation for each decoupled direction and energy group. The AP property is proved by asymptotic analysis. Numerical results of both optically thin and thick regions are presented. We can see that coarse meshes in energy, space, and time can be used to get the correct results. We will focus on the extension of the algorithm to high-dimensional problems and curvilinear geometrical problems in the future.

\bigskip
\textbf{Acknowledgement:}  X. J. Zhang and M. Tang  were
 supported by NSFC12031013, Shanghai pilot innovation project 21JC1403500 and the Strategic Priority Research Program of Chinese Academy of Sciences Grant No.XDA25010401; X. J. Zhang was supported by  China Postdoctoral Science Foundation No.2022M722107; P. Song was supported by  NSFC12031001; Y. Shi was supported by  NSFC12001052.
\bibliography{FRTE.bib}
\bibliographystyle{siam}
 \appendix
\renewcommand{\appendixname}{Appendix~\Alph{section}}
\section{The derivation of the two limiting models}\label{apen01}
In the scaling of \eqref{eqn:sca01},
we assume that the intensity $I$ and material temperature $T$ can be expanded as in \cite{MWTS}:
\begin{equation}\label{expan}
\begin{aligned}
I= I^{(0)}+\eps I^{(1)}+\eps ^2 I^{(2)}+\cdots,\\
T= T^{(0)}+\eps T^{(1)}+\eps ^2 T^{(2)}+\cdots,
\end{aligned}
\end{equation}
and the Planck function follows that
\begin{equation}\label{t4}
B= B^{(0)}+\eps B^{(1)}+\cdots,
\end{equation}
where 
$$
B^{(0)}=B|_{\eps=0}=B|_{T=T^{{(0)}}},
$$
$$
B^{(1)}=\f{\partial B}{\partial \eps}\Bigg|_{\eps=0}=\f{\partial B}{\partial T}\f{\partial T}{\partial \eps}\Bigg|_{\eps=0}=\f{\partial B}{\partial T}\Bigg|_{T=T^{{(0)}}}T^{(1)},
$$ 
By substituting  ansatzes \eqref{expan} and \eqref{t4} into  equation  \eqref{eqn:001-eq01} and collecting the terms of the same order in $\eps$, we have
\begin{subequations} \label{limit-case1}
\begin{align}
 \order \left(\f{1}{\eps}\right)&: I^{(0)}=B^{(0)},\label{limit-case11}\\
  \order \left(1\right)&: \Omega\cdot\nabla I^{(0)}=\sigma_a\left( B^{(1)} -I^{(1)}\right).\label{limit-case12}
 \end{align}
\end{subequations}
Taking \eqref{eqn:001-eq01} integral  with respect to $\nu$ from 0 to $\infty$, and taking its integral with  respect to $\Omega$, then adding the obtained equation to \eqref{eqn:001-eq02}, and dividing both sides of the equation by $\eps$, one has
$$
4\pi\partial_t\int\limits_0^{\infty}\rho \text{d}\nu+\f{1}{\eps}\int\limits_{4\pi}\int\limits_0^{\infty }\Omega \cdot \nabla I \text{d}\nu \text{d}\Omega+ C_v\partial_t T =0.
$$
By using \eqref{limit-case11}, \eqref{limit-case12} and condition $\average{\Omega\Omega}:=D_d=\f{1}{3}I_d$  (where $I_d$ denotes the 3 by 3 identity matrix), and sending $\eps\to 0$, which give
\begin{equation}\label{equi-diff}
    4\pi\partial_t\int\limits_0^{\infty}B^{(0)} \text{d}\nu+ C_v\partial_t T^{(0)} =\int\limits_0^{\infty }\nabla \cdot \left(\f{4\pi D_d}{\sigma_a}\nabla B^{(0)}\right) \text{d}\nu .
\end{equation}
by using the following relations 
$$
\int\limits_0^{\infty }4\pi B \ \text{d}\nu=T^4, \qquad  \int\limits_0^{\infty }4\pi\frac{\partial B}{\partial T} \ \text{d}\nu=4T^3,
$$
then  \eqref{equi-diff} reduces to
 \begin{equation}
 \partial_t (T^{(0)})^4+C_v \partial_t T^{(0)}=\nabla \cdot \vpran{\frac{D_d}{\sigma_R  }\nabla (T^{(0)})^4}\,,
\end{equation}
 with 
the Rosseland mean opacity $\sigma_R$ which is given by
$$
 \frac{1}{\sigma_{R}(x, T)} \equiv \frac{\int_{0}^{\infty} \frac{4\pi}{\sigma_a(x, \nu, T)} \frac{\partial B^{(0)}(\nu, T)}{\partial T} \text{d} \nu}{\int_{0}^{\infty}4\pi \frac{\partial B^{(0)}(\nu, T)}{\partial T} \text{d} \nu}=\frac{\int_{0}^{\infty} \frac{4\pi}{\sigma_a(x, \nu, T)} \frac{\partial B^{(0)}(\nu, T)}{\partial T} \text{d} \nu}{4  (T^{(0)})^{3}}.
 $$
 
In the scaling of \eqref{eqn:sca02},
the same expansions \eqref{expan} and \eqref{t4} are used, moreover, substitute \eqref{expan} and \eqref{t4} into the radiation transfer equation and collect the terms of the same order in $\eps$, which gives
\begin{subequations} \label{limit-case4}
	\begin{align}
		\order \left(\f{1}{\eps}\right)&: I^{(0)}=\rho^{(0)},\label{limit-case41}\\
		\order \left(1\right)&: \Omega\cdot\nabla I^{(0)}=\sigma_a\left( B^{(0)} -I^{(0)}\right)+\sigma_s\left( \rho^{(1)} -I^{(1)}\right).\label{limit-case42}
	\end{align}
\end{subequations}
Taking  \eqref{limit-case42}  integral  with respect to $\Omega$, one has
$$
B^{(0)}=\rho^{(0)}, 
$$
by similar calculations in the first part, one can get the following equilibrium system:
\begin{equation}
	\partial_t (T^{(0)})^4+C_v \partial_t T^{(0)}=\nabla \cdot \vpran{\frac{D_d}{\sigma_S  }\nabla (T^{(0)})^4}\,,
\end{equation}
and the mean opacity $\sigma_S$ is given by
$$
\frac{1}{\sigma_{S}} \equiv \frac{\int_{0}^{\infty} \frac{4\pi}{\sigma_s} \frac{\partial B^{(0)}(\nu, T)}{\partial T} \text{d} \nu}{\int_{0}^{\infty}4\pi \frac{\partial B^{(0)}(\nu, T)}{\partial T} \text{d} \nu}=\frac{\int_{0}^{\infty} \frac{4\pi}{\sigma_s} \frac{\partial B^{(0)}(\nu, T)}{\partial T} \text{d} \nu}{4  (T^{(0)})^{3}}.
$$ 

In the scaling of \eqref{eqn:sca03},
the same expansions \eqref{expan} and \eqref{t4} are used, moreover, substitute \eqref{expan} and \eqref{t4} into the radiation transfer equation and collect the terms of the same order in $\eps$, which gives
\begin{subequations} \label{limit-case3}
\begin{align}
 \order \left(\f{1}{\eps}\right)&: (\sigma_a+\sigma_s)I^{(0)}=\sigma_a B^{(0)}+\sigma_s\rho^{(0)},\label{limit-case31}\\
  \order \left(1\right)&: \Omega\cdot\nabla I^{(0)}=\sigma_a\left( B^{(1)} -I^{(1)}\right)+\sigma_s\left( \rho^{(1)} -I^{(1)}\right).\label{limit-case32}
 \end{align}
\end{subequations}
Taking \eqref{limit-case31} and \eqref{limit-case32}  integral  with respect to $\Omega$, one has
$$
B^{(0)}=\rho^{(0)}, \qquad B^{(1)}=\rho^{(1)},
$$
by similar calculations in the first part, one can get the following equilibrium system:
\begin{equation}
 \partial_t (T^{(0)})^4+C_v \partial_t T^{(0)}=\nabla \cdot \vpran{\frac{D_d}{\sigma_T  }\nabla (T^{(0)})^4}\,,
\end{equation}
and the mean opacity $\sigma_T$ is given by
$$
 \frac{1}{\sigma_{T}} \equiv \frac{\int_{0}^{\infty} \frac{4\pi}{(\sigma_a+\sigma_s)} \frac{\partial B^{(0)}(\nu, T)}{\partial T} \text{d} \nu}{\int_{0}^{\infty}4\pi \frac{\partial B^{(0)}(\nu, T)}{\partial T} \text{d} \nu}=\frac{\int_{0}^{\infty} \frac{4\pi}{(\sigma_a+\sigma_s)} \frac{\partial B^{(0)}(\nu, T)}{\partial T} \text{d} \nu}{4  (T^{(0)})^{3}}.
 $$

 In the frequency-dependent radiative diffusion regime, the same expansions \eqref{expan} and \eqref{t4} are used, moreover, substitute \eqref{expan} and \eqref{t4} into equation \eqref{eqn:001-noneq01} and collect the terms of the same order in $\eps$, which gives
 \begin{subequations} \label{limit-case2}
 	\begin{align}
 	\order \left(\f{1}{\eps}\right)&: I^{(0)}=\rho^{(0)},\label{limit-case21}\\
 	\order \left(1\right)&: \Omega\cdot\nabla I^{(0)}=\sigma_s\left( \rho^{(1)} -I^{(1)}\right).\label{limit-case22}
 	\end{align}
 \end{subequations}
 Taking \eqref{eqn:001-noneq01}  integral  with respect to $\Omega$,  and dividing both sides of the equation by $\eps$, one has
 $$
 4\pi\partial_t\rho +\f{4\pi}{\eps}\average{\Omega \cdot \nabla I}=4\pi\sigma_a(B-\rho),
 $$ 
 by similar calculations in the first part, one can get the following non-equilibrium system:
 \begin{equation}
 \left\{
 \begin{aligned}
 &\partial_t \rho^{(0)}- \nabla \cdot \vpran{\frac{D_d}{\sigma_{s}  }\nabla \rho^{(0)}}=\sigma_a(B^{(0)}-\rho^{(0)}),\\
 &  C_v\partial_t T^{(0)}  =\int_{0}^{\infty}4\pi\sigma_a(\rho^{(0)}-B^{(0)}) \text{d} \nu. \ \  \\
 \end{aligned}
 \right.
 \end{equation} 
\end{document}